\documentclass[11pt]{elsarticle}

\usepackage{upgreek}
\usepackage{libertine} 
\usepackage{euler} 
\usepackage{eucal} 
\usepackage[utf8]{inputenc}

\usepackage[T1]{fontenc}
\usepackage{textcomp}
\usepackage[utf8]{inputenc}

\usepackage{amsopn, amsthm, amsgen, amscd,amsmath,amssymb}
\usepackage[bbgreekl]{mathbbol}
\DeclareMathAlphabet{\mathbbm}{U}{bbm}{m}{n}
\usepackage{color}
\usepackage{tikz}
\usepackage{epstopdf}

\usepackage{graphicx,pspicture}

\newcommand{\ojo}[1]{%
   \hskip 1sp \marginpar{\small\sffamily\raggedright
     SAH: #1}}
\newcommand{\ojoA}[1]{%
   \hskip 1sp \marginpar{\small\tt\raggedright
     Andr\'es: #1}}

\define\cf{\; {\rm cf} \; }

\define\pmin{{\cal P}^-}

\define\Dom{\; {\rm Dom} \; }
\define\ran{\; {\rm rng} \; }

\define\Mod{\; {\rm Mod} \; }

\def\cK{{\mathcal{K}}}

\def\b{\beta}
\def\a{\alpha}
\def\l{\lambda}

\def\z{\zeta}

\def\w{\omega}

\def\g{\gamma}

\def\al{\aleph}

\font\frac=eurb10

\def\k{\kappa}

\def\gd{{\frak D}}

\define\qed{\hfill $\Box$}

\newcommand{\smallbox}[1]{\leavevmode\thinspace\hbox{\vrule\vtop{\vbox
   {\hrule\kern1pt\hbox{\strut\thinspace{#1}\thinspace}}
   \kern1pt\hrule}\vrule}\thinspace}

\newcount\skewfactor
\def\mathunderaccent#1#2 {\let\theaccent#1\skewfactor#2
\mathpalette\putaccentunder}
\def\putaccentunder#1#2{\oalign{$#1#2$\crcr\hidewidth
\vbox to.2ex{\hbox{$#1\skew\skewfactor\theaccent{}$}\vss}\hidewidth}}

\newcount\skewfactor
\def\mathunderaccent#1#2 {\let\theaccent#1\skewfactor#2
\mathpalette\putaccentunder}
\def\putaccentunder#1#2{\oalign{$#1#2$\crcr\hidewidth
\vbox to.2ex{\hbox{$#1\skew\skewfactor\theaccent{}$}\vss}\hidewidth}}

\makeatletter
\def\indsym#1#2{%
  \setbox0=\hbox{$\m@th#1x$}%
  \kern\wd0%
  \hbox to 0pt{\hss$\m@th#1\mid$\hbox to 0pt{$\m@th#1^{#2}$}\hss}%
  \lower.9\ht0\hbox to 0pt{\hss$\m@th#1\smile$\hss}%
  \kern\wd0} 
\def\nindsym#1#2{%
  \setbox0=\hbox{$\m@th#1x$}%
  \kern\wd0%
  \hbox to 0pt{\mathchardef\nn="3236\hss$\m@th#1\nn$\kern1.4\wd0\hss}
  \hbox to 0pt{\hss$\m@th#1\mid$\hbox to 0pt{$\m@th#1^{#2}$}\hss}%
  \lower.9\ht0\hbox to 0pt{\hss$\m@th#1\smile$\hss}%
  \kern\wd0}

\makeatother

\newcommand{\chxyz}{(\bar{x},\bar{y},\bar{z})}

\newcommand{\bbz}{{\Bbb Z}}

\newcommand{\ialak}{[I]^k}
\newcommand{\ialakmasuno}{[I]^{k+1}}

\newcommand{\Proof}{\noindent{\sc Proof} \hspace{0.2in}}

\newcommand{\lesdot}{\mathrel{\mathord{<}\!\!\raise 0.8
pt\hbox{$\scriptstyle\circ$}}}

\newcommand{\GK}{{\frak K}}

\newtheorem{theorem}{Theorem}[section] 
\newtheorem{lemma}[theorem]{Lemma} 
\newtheorem{proposition}[theorem]{Proposition} 
\newtheorem{corollary}[theorem]{Corollary} 
 
\newtheorem{definition}[theorem]{Definition}

\newtheorem{claim}[theorem]{Claim}
\newtheorem{conclusion}[theorem]{Conclusion}

\newtheorem{remark}[theorem]{Remark} 
 
\newtheorem{fact}[theorem]{Fact} 
 
\newtheorem{discussion}[theorem]{Discussion}
\newtheorem{hypothesis}[theorem]{Hypothesis}

\newtheorem{context}[theorem]{Context}

\begin{document}

\begin{frontmatter}

	\title{The Hart-Shelah example, in stronger logics.}
\author[rutg,jer]{Saharon~Shelah\fnref{fn1}}
\ead{shelah@math.huji.ac.il}
\author[unal]{Andrés~Villaveces\fnref{cor1,fn2}}
\ead{avillavecesn@unal.edu.co}
\fntext[cor1]{Corresponding author}
\fntext[fn1]{The first author's research was partially supported
by `BSF' (USA-Israel); publication no 648 in the first author's
publication list.}
\fntext[fn2]{The second author was sponsored by
  Colciencias. Part of the research
  was done during a visit to the Centre de Recerca Matemàtica
  (Bellaterra, Catalonia)'s Intensive Research Program in Strong
  Logics and Large Cardinals, 2016.}
  \address[rutg]{Department of Mathematics, Rutgers University,
	  Hill Center for the Mathematical Sciences, 110 Frelinghuysen Rd.,
  Piscataway, NJ 08854-8019}
  \address[jer]{Einstein Institute of Mathematics,
	  Edmond J. Safra Campus,
The Hebrew University of Jerusalem,
Givat Ram.
Jerusalem, 9190401, Israel}
\address[unal]{Departamento de Matemáticas, Universidad Nacional de Colombia,
	Av. Carrera 30 \# 45-03. Ciudad Universitaria.
Bogotá, Colombia - 111321}

\begin{abstract}
We generalize the Hart-Shelah example~\cite{HaSh:323} to
higher infinitary logics. We build, for each natural number $k\geq 2$ and for each
infinite cardinal $\lambda$, a sentence $\psi_k^\lambda$ of the logic
$L_{(2^\lambda)^+,\omega}$ that (modulo mild set theoretical
hypotheses around $\lambda$ and assuming $2^\l < \l^{+m}$) is
categorical in
$\lambda^+,\dots,\lambda^{+k-1}$ but not in
$\beth_{k+1}(\lambda)^+$ (or beyond); we study the dimensional encoding of
combinatorics involved in the construction of this sentence and study
various model-theoretic properties of the resulting abstract
elementary class
$\cK^*(\lambda,k)=(Mod(\psi_k^\lambda),\prec_{(2^\lambda)^+,\omega})$
in the finite interval of cardinals
$\lambda,\lambda^+,\dots,\lambda^{+k}$.
\end{abstract}


\end{frontmatter}



The study of categoricity transfer has been central to model theory since
Morley's theorem; the question of finding extensions of this theorem to
infinitary contexts and to abstract elementary classes has been a major source
of results. Many central concepts of stability theory, both in first order and
in its generalizations, are essential byproducts of the theory built in order
to generalize the original Morley theorem.


One of the most important landmarks along this path was the Categoricity
Transfer result for $L_{\omega_1,\omega}$ due to the first author: if a
sentence $\psi$ is categorical in $\aleph_n$ for all $n<\omega$ and the weak
GCH holds for the $\aleph_n$'s ($2^{\aleph_n}<2^{\aleph_{n+1}}$ for all
$n<\omega$) then $\psi$ is categorical in all cardinals
(see~\cite{Sh:87a} and~\cite{Sh:87b}; although these are two references, they
correspond to ``Part A'' and ``Part B'' of one big paper from 1983).
Notice the unusually strong assumption!

An example from 1990 due to Bradd Hart and the first author of this
paper~\cite{HaSh:323} established the
(surprising) necessity of that strong assumption: the existence of \emph{few models} at all the
$\aleph_n$'s is needed to get the eventual categoricity transfer for
$L_{\omega_1,\omega}$: they provide, for each positive $k\in \omega$,
an example of a
sentence $\psi_k$ in $L_{\omega_1,\omega}$  categorical in
$\aleph_0,\aleph_1,\cdots,\aleph_k$ but not eventually categorical:
there exists some cardinal greater than $\aleph_k$ where categoricity fails.

That important example has later been
referred to as \textbf{the Hart-Shelah example.} In many ways, the existence of such sentences
points to an interesting failure of the ``categoricity transfer'' (Morley's theorem for first
order logic, for countable theories) at \emph{small} cardinalities, in the absence of further set
theoretical hypotheses. Our work extends those results.

The present paper
shows that a similar example also exists for the stronger logic
$L_{(2^\lambda)^+,\omega}$. A corollary of our result is that
any extension of the results from~\cite{Sh:87a} and~\cite{Sh:87b} 
to stronger logics will require to assume categoricity at \textbf{all} cardinalities
$\lambda,\lambda^+,\dots,\lambda^{+n},\dots$ for all $n<\omega$.

Later, the first author has attempted an extension of the main result
from~\cite{Sh:87a} and~\cite{Sh:87b} to Abstract Elementary Classes. These are
more general than classes axiomatized by the logic $L_{(2^\lambda)^+,\omega}$. 

%

Our construction provides, for each infinite cardinal $\lambda$ and each $k\in
(2,\omega)$, a sentence $\psi_k$ of
$L_{(2^\lambda)^+,\omega}$ that is categorical in
$\lambda,\lambda^+,\dots,\lambda^{+k}$ but is not categorical in any
cardinality $\mu\geq \beth_{k+1}(\lambda)^+$.



The shift of focus from infinitary logic to abstract elementary
classes entails in many cases using Galois (orbital) types instead of
syntactic types; although this shift is natural, compactness and
locality properties in general do not transfer to Galois types. In
particular, \emph{tameness} and \emph{type-shortness} do not hold in
general for Galois types. Tameness was isolated by Grossberg and
VanDieren~\cite{GrVa}; later, Baldwin and Shelah~\cite{BlSh} constructed an
example of failure of tameness, based on an almost free non-Whitehead group.
More recently, Boney and Unger have provided serious set theoretic reasons for
the failure of tameness in AECs~\cite{BonUng17}.

In~\cite{BlKo}, Baldwin and Kolesnikov study again the Hart-Shelah
example: they prove that for the sentence $\psi_k$ of
$L_{\omega_1,\omega}$ of the example, the corresponding AEC (for
$k\geq 3$)
\[\cK^{HS}(\omega_1,k)=(Mod(\psi_k),\prec_{\omega_1,\omega})\]
\begin{itemize}
\item has disjoint amalgamation,
\item is Galois stable exactly in
  $\aleph_0,\aleph_1,\dots,\aleph_{k-1}$,
\item is $(<\aleph_0,\leq \aleph_{k-1})$-tame.
\end{itemize}
Moreover, the AEC axiomatized by their sentence $\psi_k$ fails 
$(\aleph_{k-1},\aleph_k)$-tameness. This is an immediate consequence of the
failure of categoricity transfer and the upward categoricity theorem for tame
AECs due to Grossberg and VanDieren~\cite{GrVa2}.

Baldwin and Kolesnikov really study a slight variant
of the Hart-Shelah example, presented in the language of group actions
and revealing the filiation to the early Baldwin-Lachlan example of an
$\aleph_1$-categorical theory which is not almost strongly
minimal.

More recently, Boney~\cite{Bon14a} has continued this study of
the behavior of the Hart-Shelah example; he has proved that the class
$\cK^{HS}(\omega_1,k)$ has a ``good $\aleph_m$-frame'' for all $m\leq
k-1$ but cannot have a good frame above by the failure of
stability. Then, Boney and Vasey~\cite{BoVa} continue this study and
show first that the frame at $\aleph_{k-1}$ cannot be
``successful''. They study good frames in connection with the
Hart-Shelah example: for frames around the $\aleph_n$'s ($n<\omega$)
the Hart-Shelah example is a natural place to look for ``boundary
properties'': being ``successful up to some point'' but failing to be
successful above.

Our generalization of the Hart-Shelah example addresses the question of how
necessary an assumption similar to ``few models in all the $\aleph_n$'s'' is
for categoricity transfer in the case of stronger logics. Here of course the
corresponding assumption would be of the form ``few models in all the
$\lambda^{+n}$ ($n<\omega$)''.


We build a sentence $\psi_k^\lambda$ in $L_{(2^\lambda)^+,\omega}$,
categorical in $\lambda,\lambda^+,\dots,\lambda^{+k}$ but not
categorical in $\beth_{k+1}(\lambda)^+$. Here are two important differences
between our approach and earlier ones:
\begin{itemize}
\item The sentences are constructed in all cases by first building a
  ``standard model'' and then extracting the sequence from it. In the
  Hart-Shelah example, one predicate $Q$ ``ties together'' various
  copies of groups in a way that ends up linking the ``dimension'' of
  the predicate to the length of induction in the proof of
  categoricity. In our example, we need a large family of predicates
  $Q_s$, $s\in S=[\lambda]^{<\aleph_0}$.
\item The ``failure of categoricity'' argument at cardinals greater
  than or equal to $\beth_{k+1}(\lambda)^+$ here is done by using a
  regular filter $\mathfrak D$.
\end{itemize}

A natural question arises, on the ``gap'' between categoricity and failure of
categoricity of $\psi_k^\lambda$. Here, we can guarantee categoricity in the
interval $[\lambda,\lambda^{+k}]$ and failure of categoricity\dots at
$\beth_{k+1}(\lambda)^+$. Admittedly, this is a very large gap, relatively much
wider than what Baldwin and Kolesnikov have for their version of the
Hart-Shelah sentence. The question remains open whether this gap may be
reduced.

In our concluding remarks, we raise some questions connected with the tameness
and frames, inspired by the paper~\cite{BoVa}. In particular, we ask whether
the methods from that paper (that worked for the Hart-Shelah sentence) may be
generalized to our sentence $\psi_k^\lambda$.

%
%


\underline{A note on indexing:} the previous papers dealing with constructing
examples of sentences where categoricity ``stops''
are~\cite{HaSh:323},~\cite{BlKo} (which
proved more model theoretic facts on a variant of the original example
and studied the abstract elementary class determined by the example;
in particular, Galois (=orbital) types, the amalgamation and tameness
spectra associated with the class),~\cite{Bon14a} and~\cite{BoVa}, in
which the connection to frames is worked out (analyzing the
Hart-Shelah example enables Boney and Vasey to study limitations to
the existence of good frames). Now, for~\cite{HaSh:323}, the
``critical'' cardinality (the last cardinality of categoricity) is
$\aleph_k$. In~\cite{BlKo}, because of the way they analyze the
construction, it is more natural to work with $k\geq 3$ and with $k-2$
as the critical cardinality. The two other papers follow this.

Since our paper is directly a generalization of~\cite{HaSh:323}, it is
more natural for us to revert to the choice of critical cardinality
from there, of course adapted to our context. So, the last cardinality
where we will have categoricity is $\lambda^{+k-1}$.

Our notation is standard.

We thank John Baldwin, Will Boney, Rami Grossberg, Alexei Kolesnikov,
Sebastien Vasey and Boban
Velickovic for several remarks and valuable discussions
concerning (directly or less directly) this work, as well as for
pressing us to provide some clarification of the big
construction. The second author is particularly indebted to John
Baldwin for very interesting conversations of the connections between
this example and the original group covers in the strongly minimal
context, due to Baldwin and Lachlan~\cite{BlLa71}. We also thank Péter Komjath for helpful
discussion on the negative partition relation in~\cite{EHMR} useful in our theorem.
We also thank the anonymous referee of an earlier version of this paper for
extremely insightful and helpful comments. They (hopefully) led to our
improving this paper. We also thank a second anonymous referee of the version prior to this for remarks that
led to a substantial rewriting and what we believe is a much better presentation of the results.


\section{Construction of the sentence $\psi^\lambda_k$, in $L_{(2^\lambda)^+,\omega}$}

\begin{context}
    For the rest of the paper, we fix an infinite cardinal $\lambda$ and a natural number
$k\geq 2$.  
\end{context}

We build in this section a new sentence $\psi^\lambda_k$ in the logic $L_{(2^\lambda)^+,\omega}$.
Our construction of $\psi^\lambda_k$ requires first building a model we will call ``canonical'',
$M_I$, for an arbitrary index set $I$ and later taking a conjunction of the first order theory of
$M_I$ along with several infinitary sentences describing the behavior of various components of
$M_I$. The sentence $\psi^\lambda_k$ has some similarity to the Hart-Shelah sentence and may be seen
as a generalization, but important differences are also present and will be apparent later
(the regular group $G$ and the regular filter on $\lambda$, $\gd$). However, it is important
to stress that prior knowledge of the Hart-Shelah is \emph{not necessary} for an understanding of
our construction, as we make it self-contained.

We will build the model $M_I$ around a ``spine'' $I$, essentially by coding interactions between
$k$-element subsets of $(k+1)$-element subsets of $I$, in some combinatorial ways. Namely, we will
define various groups and encode in the model \emph{their actions} on those $k$ and $(k+1)$-element
subsets of $I$, focusing especially on the way different $k$-subsets of a given $(k+1)$-subset of $I$
interact. Finally, a collection of predicates (called $Q_s$) will ``tie'' those combinatorial
interactions.

\begin{definition}\label{bricks}{\bf Notation and general construction tools.}
We fix the following basic objects to use in the construction later.
\begin{itemize}
\item $S = S_\l := [\l]^{<\al_0} = \{u\subset \l | u\mbox{ is finite}\}$,
\item $\gd = \gd_\l := \{ A\subset S | \exists u_A\in S \quad \forall v\in
S(u_A\subset v \to
v\in A)\}$, the \emph{regular filter} on $S$ generated by sets of the form
$\langle u\rangle = \{ v\in S | u\subset v\}$,
\item $G^+ = G^+_\l := {}^S(\bbz_2)$, as a group with the natural operation
$(f+g)(v) = f(v)+_{\bbz_2}g(v)$,
\item $G = G_\l := \{ f\in {}^S(\bbz_2) | \ker(f) = \{ u\in S | f(u) =
0 \} \in \gd
\}$, as a subgroup of $G^+$ ($G\leq G^+$ since, if $f,g \in G$, then
$\ker(f),\ker(g)\in \gd$, so $\ker(f+g) \supset \ker(f)\cap \ker(g)\in
\gd$, hence $\ker(f+g)\in \gd$ and $f+g\in G$.
\end{itemize}
\end{definition}

Note that $|G|=2^\lambda$.

\subsection{The model $M_I$}\label{spine}

Fix a set $I$ for the rest of this section.

We define the group $H_I$ and then describe the model $M_I$.

\begin{definition}[The group $H_I$]\label{groupHI}
    For our (fixed) $I$, we let $H_I := [\ialak]^{<\al_0}$. So, $H_I$ is the set of finite
    subsets of $\ialak$, with the group
	operation $F_1+F_2:= F_1\Delta F_2$ (symmetric difference).
    Equivalently, $H_I$ may be
	seen as the set of functions $\ialak\to\bbz_2$ with finite support. In this case, $F\in
	H_I$ is coded by the function $h_F:\ialak\to\bbz_2$ with $h_F(u)=1$ iff $u\in F$ and the
	group operation is given by $(h_1+h_2)(u)=h_1(u)+_{\bbz_2}h_2(u)$.
\end{definition}
    
    We now start the (lengthy) description of $M_I$: the universe, basic predicates, projections
    between them, other partial functions coded by relations and, crucially, \textbf{the family of
    predicates} $Q_s$ (for $s\in S$).

    \begin{definition}[Universe of $M_I$]
	The universe of $M_I$ is the union of seven different sorts:
\[ |M_I| = I \cup \ \ialak \cup \ \ialakmasuno \cup (\ialak \times
S\times H_I )\cup (\ialak \times S\times \bbz_2) \cup H_I \cup (\ialakmasuno
\times G).\] 
	\label{universeM_I}
    \end{definition}

The following remarks on the universe of $M_I$ are important:

\begin{itemize}
\item The natural way to think about the universe of $M_I$ is as
\[ |M_I|\]
\begin{center}
    consisting of two parts: the ``support of the model'' ($I$, $\ialak$, $\ialakmasuno$) and many
    copies of (the domains of) the three groups $H_I$, $\bbz_2$ and $G$, \emph{indexed} by elements of
    $S$ and of the support part: 

\end{center}
\[ \overbrace{I \cup \ \ialak \cup \
\ialakmasuno}^{\mbox{\footnotesize `support part'}}
 \cup \overbrace{(\ialak \times
S\times H_I )\cup (\ialak \times S\times \bbz_2) \cup H_I \cup (\ialakmasuno
\times G)}^{\mbox{\footnotesize $H_I$ and copies of the domains of $H_I$, $\bbz_2$
and $G$}}.\] 
\item Notice that the intersection between all the pieces of the
  model is empty.
\item The universe of $M_I$ depends directly on $I$ and on $G$, as is
  clear from
  the various pieces. In particular, when the cardinality of $I$ is
  $\geq\lambda$, the cardinality of $M_I$ will be equal to
  $|I|+2^\lambda$.
\item The universe depends on $k$ as well. Of
  course in our standard model this dependence is immediate, as seen
  from the superindices $k$ and $k+1$. In general models later, we will need projection
  functions among the
  predicates in the model in order to axiomatize the connections
  between pieces corresponding to abstract versions of $I$, $[I]^k$, etc. \emph{This dependence on
    $k$ will be crucial in the ``dimension'' analysis later.}
\item The universe also depends on $\lambda$, through the
  appearance of $S$ and $G$ among the pieces.
\end{itemize}

\begin{definition}[Relations, functions of $M_I$ - the predicates $Q_s$]
    The structure of $M_I$ consists of the following items:
    \begin{itemize}
	\item $\lambda$-many predicates $P_0^M,P_{1,1}^M,P_{1,2}^M,P_2^M,(P_{2,s}^M)_{s\in
	    S},P_3^M,(P_{3,s}^M)_{s\in S}, P_4^M,P_5^M$,
	\item $k$-many projections $\pi^0_\ell:P_{1,1}^M\to P_0^M$ ($\ell<k$) and $k+1$-many
	    projections $\pi_\ell^1:P^M_{1,2}\to P_0^M$,
	\item $2^\lambda$-many additional functions $F_2^M,F_3^M,F_4^M,F_5^M,(F^M_{3,g^*})_{g^*\in
	    G}$,
	\item a $(3k+4)$-ary predicate $Q_s$, for each $s\in S$.
    \end{itemize}
    Each of these predicates and functions will be discussed in detail in the following paragraphs.
    \label{structureM_I}
\end{definition}

\subsubsection{Descriptions of basic relations, functions, and the $Q_s$-predicates}

\begin{description}
\item[Basic Relations:] these consist of a family of $\lambda$-many
  predicates
\[P_0^M,P_{1,1}^M,P_{1,2}^M,P_2^M,(P_{2,s}^M)_{s\in S}, 
P_3^M,(P_{3,s}^M)_{s\in S}, P_4^M,P_5^M\]
defined by
\begin{itemize}
\item $P_0^M = I$,
\item $P_{1,1}^M = \ialak$,
\item $P_{1,2}^M = \ialakmasuno$,
\item $P_2^M = \ialak \times S\times H_I$,
\item for $s\in S$, $P_{2,s}^M = \{ (u,s,h) \in P_2^M | u\in \ialak,
h\in H_I\} = \ialak \times \{ s\} \times H_I$,
\item $P_3^M = \ialak \times S\times \bbz_2$ (a copy of $\bbz_2$ for each
$b\in \ialak$, $s\in S$),
\item for $s\in S$, $P_{3,s}^M = \{ (u,s,i) \in P_3^M | u\in \ialak,
i\in \bbz_2\} = \ialak \times \{ s\} \times \bbz_2$,
\item $P_4^M = H_I$,
\item $P_5^M = \ialakmasuno\times G$
\end{itemize}
\begin{quotation}

  \begin{remark}
The meaning of $P_0^M,P_{1,1}^M,P_{1,2}^M,P_2^M$, $P_3^M$, $P_4^M$ is
clear. In the case of $P_{2,s}^M$, the idea is that we stack
``copies'' of $H_I$ for each $b\in \ialak$ and each $s\in S$, and
similarly for $P_3^M$, $P_{3,s}^M$. Another way of seeing this is
thinking of the predicates as codifying families, as follows:
\begin{itemize}
\item $P_2^M$ corresponds to $(H_{v,s})_{v\in \ialak,s\in S}$, 
\item $P_3^M$ corresponds to $((\bbz_2)_{v,s})_{v\in \ialak,s\in S}$,
\item $P_5^M$ corresponds to $(G_u)_{u\in \ialakmasuno}$.
\end{itemize}    
  \end{remark}

\end{quotation}
\item[Projections:] 
We also include, for $\ell <k$, all the projections
$\pi^0_\ell:P_{1,1}^M \to P_0^M$:
\[\pi^0_\ell(\bar{a}) = a_\ell,\]
and for $\ell < k+1$, the projections
$\pi^1_\ell:P_{1,2}^M \to P_0^M$:
\[\pi^1_\ell(\bar{a}) = a_\ell.\]
The role of these projections is to tie the predicates 
$P_{1,1}^M$ and $P_{1,2}^M$ to $P_0^M$ making them behave as the
corresponding sets of $k$tuples or $k+1$-tuples.
\item[Other Partial Functions:]
We also include $2^\lambda$-many functions in $M_I$,
\[F_2^M,F_3^M,F_4^M,F_5^M,(F^M_{3,g^*})_{g^*\in G}:\]
\begin{itemize}
\item A unary function $F_2^M$ with domain $P_2^M$, given by
\[F_2^M(u,s,h) = u,\]
\item A unary function $F_3^M$ with domain $P_3^M$, given by
\[F_3^M(u,s,i) = u,\]
\item for $g^*\in G$, a unary function $F^M_{3,g^*}$ with domain
$P^M_5$, given by \[F^M_{3,g^*}(u,g) = (u,g^*+g),\]
\item A binary function $F^M_4$ with domain $P^M_2\times
P_4^M$, given by 
\[F_4^M\Big( (v,s,h),h_1\Big) = (v,s,h+_Hh_1)\]
\item A unary function $F^M_5$ with domain $P^M_5$, given by 
\[F_5^M(u,g) = u,\]
\end{itemize}





\item[A $(3k+4)$-ary predicate] $Q_s$, for each $s\in S$. This is the
  crux of the construction of the model $M_I$. The predicate will
  encode interactions between the different parts of the model, in a
  way that will involve \emph{dimensional} interactions between
  them. This predicate on the one hand \emph{enables} later to move up
  in the proof of categoricity by induction $k-1$ times from $\lambda$
  to $\lambda^k$ and on the other \emph{blocks} the proof from moving
  up to $\lambda^{k+1}$.
It is interpreted in $M_I$ as
the set of tuples\label{defM_I} 
\[ \langle a_0,\dots,a_k,u_0,\dots,u_k,x_0,\dots,x_{k-1},y_k,z \rangle \]
satisfying (for fixed $s\in S$!!) for all $h_k \in H_I$, $i_\ell \in
\bbz_2 (\ell < k)$, $g\in G$:
\begin{description}
\item[$(\a)$] $a_\ell \in I$ with no repetitions ($\ell \leq k$),
\item[$(\b)$] $u_\ell=\langle a_m |m\not= \ell\rangle\in P_{1,1}^M$
($\ell \leq k$),
\item[$(\g)$] $y_k = (u_k,s,h_k)\in P_2^M$,
\item[$(\delta)$] $x_\ell$ has the form $(u_\ell,s,i_\ell)\in P_3^M$
($\ell < k$) so $i_\ell\in\bbz_2$, 
\item[$(\epsilon)$] $z$ is of the form $(u,g)\in P_5^M$,
where $u=(a_0,\dots,a_k)\in \ialakmasuno$ and
\item[$(\z)$] {\bf (main point)}
\[\bbz _2 \models \sum_{\ell<k}i_\ell = h_k(u_0) + g(s).\]
\end{description}
\end{description}

Some general remarks on this definition of the model $M_I$ are in point, before giving a specific
description for the case $k=2$.

\begin{remark}
  \begin{itemize}
      \item $(\z)$ is the crucial part of the definition of the predicates $Q_s$. It provides the
	  connection between $k$ copies of $\bbz_2$, one copy of $H_I$, one copy of $G$ and the
	  $(k+1)$-many $k$-element subsets of a set of size $k+1$ in $I$.
      \item The role of $F_2^M$ is to project $P_2^M$ (essentially
    $(H_{v,s})_{v\in \ialak,s\in S}$) onto its first coordinate; to
    trace the $k$-element subset of $I$ it corresponds to. Similarly
    for $F_3^M$ and $F_5^M$.
\item The functions $F_{3,g^*}^M$ and $F_4^M$ encode the \emph{actions} of the groups $G$ and $H_I$
    on the corresponding ``fibers'' over $u\in [I]^{k+1}$ or $(v,s)\in [I]^k\times S$. The model
    $M_I$ does not really include the group operations corresponding to $G$ and $H_I$; it only has
    the effect of the group actions on the appropriate fibers.
  \item Notice that $+^H$ is definable - so in this case there is no
    need to add an analogue of $F_4$ for copies of $\bbz_2$:
    \[F_4^M(F_4^M ((u,s,h),h_1),h_2) =
      F_4^M((u,s,h),h_3)\Leftrightarrow H\models h_1+h_2=h_3.\] 
  \end{itemize}
\end{remark}

\subsubsection{Illustration of the definition of $M_I$, when $k=2$}

As an example to visualize the situation, we momentarily fix $k=2$. We
also fix $s\in S$ and choose some $u\in \ialakmasuno=[I]^{2+1}$,
$u=\langle a_0,a_1,a_2\rangle$. This determines automatically (using
the projections) in the models we have described so far
$a_0,a_1,a_2$ and $u_0=\langle a_1,a_2\rangle$,
$u_1=\langle a_0,a_2\rangle$, $u_2=\langle a_0,a_1\rangle$.

We then have
\begin{itemize}
    \item copies of $\bbz_2$ over both $u_0$ and $u_1$,
    \item a copy of the domain of $H_I$ over $u_2$, together with the action of $H_I$ on this copy,
    \item a copy of $G$ over $u$, again with the action of $G$ over this copy.
\end{itemize}

Furthermore, we have the predicate $Q_s$: it is in this case $3\cdot 2+4=10$-ary. The $10$-uple
associated with our $u$ is then of the form
\[ (a_0,a_1,a_2,u_0,u_1,u_2,x_0,x_1,y_2,z), \]

\noindent
with $x_0=(u_0,s,i_0)$, $x_1=(u_1,s,i_1)$, $y=(u_2,s,h_2)$ and $z=(u,g)$ for some $i_0,i_1\in
\bbz_2$, $h_2\in H_I$ and $g\in G$.

We want to describe when this tuple belongs to $Q_s$. The following triangle summarizes the relevant
information:

\begin{center}

\tikzset{every picture/.style={line width=0.75pt}} 

\begin{tikzpicture}[x=0.75pt,y=0.75pt,yscale=-.67,xscale=.67]

\draw   (323.75,148) -- (442.5,332) -- (205,332) -- cycle ;
\draw  [fill={rgb, 255:red, 0; green, 0; blue, 0 }  ,fill opacity=1 ] (317.5,148) .. controls (317.5,144.55) and (320.3,141.75) .. (323.75,141.75) .. controls (327.2,141.75) and (330,144.55) .. (330,148) .. controls (330,151.45) and (327.2,154.25) .. (323.75,154.25) .. controls (320.3,154.25) and (317.5,151.45) .. (317.5,148) -- cycle ;
\draw  [fill={rgb, 255:red, 0; green, 0; blue, 0 }  ,fill opacity=1 ] (198.75,332) .. controls (198.75,328.55) and (201.55,325.75) .. (205,325.75) .. controls (208.45,325.75) and (211.25,328.55) .. (211.25,332) .. controls (211.25,335.45) and (208.45,338.25) .. (205,338.25) .. controls (201.55,338.25) and (198.75,335.45) .. (198.75,332) -- cycle ;
\draw  [fill={rgb, 255:red, 0; green, 0; blue, 0 }  ,fill opacity=1 ] (436.25,332) .. controls (436.25,328.55) and (439.05,325.75) .. (442.5,325.75) .. controls (445.95,325.75) and (448.75,328.55) .. (448.75,332) .. controls (448.75,335.45) and (445.95,338.25) .. (442.5,338.25) .. controls (439.05,338.25) and (436.25,335.45) .. (436.25,332) -- cycle ;

\draw (341.33,273) node  [align=left] {$g$};
\draw (323,236.33) node [scale=1.44] [align=left] {$u$};
\draw (244,222.67) node [rotate=-305.13] [scale=.85]  [align=left] {$h_2:[I]^2 \to \bbz_2$};
\draw (409,235) node [rotate=-55.76] [align=left] {$i_1 \, (0/1)$};
\draw (326,348.33) node  [align=left] {$i_0 \, (0/1)$};
\draw (283,238) node [rotate=-305.13] [align=left] {$u_2$};
\draw (373,247) node [rotate=-55.76] [align=left] {$u_1$};
\draw (322,319.67) node  [align=left] {$u_0$};
\draw (466,348) node  [align=left] {$a_2$};
\draw (177,341) node  [align=left] {$a_1$};
\draw (327,126) node  [align=left] {$a_0$};

\end{tikzpicture}

\end{center}

The tuple $(a_0,a_1,a_2,u_0,u_1,u_2,x_0,x_1,y_2,z)$ belongs to $Q_s$ if and only if
\[ \bbz_2\models i_0+i_1=h_2(u_0)+g(s).\]
Therefore, on top of the triangle $u$ we have (when $k=2$) four pieces of information playing: two
elements ($i_0,i_1$) of $\bbz_2$ associated to two sides of the triangle, one element $h$ of $H_I$
associated to the third side of the triangle (and the value of $h$ at $u_0$) and finally one element
$g$ of $G$ associated to the triangle $u$ itself - and the value of $g$ at\dots $s$.

\subsection{The language, the sentence $\psi_k^\lambda$ and the AEC
  $\cK^*(\lambda,k)$}
\label{sec:languageetal}

We now build the sentence $\psi_k^\lambda$.

\begin{definition}\label{la6.2}
We deal with two vocabularies:
  \begin{itemize}
  \item Let $\tau^-$ be the vocabulary of all the construction above,
    except the predicates $\{ Q_s|s\in S\}$ and
  \item let $\tau$ be the full vocabulary used in the construction of $M_I$.
  \end{itemize}
 Specifically,

\begin{description}
\item let $\tau^- = \langle P_0,P_{1,1},P_{1,2},P_2,(P_{2,s})_{s\in
S},P_3,(P_{3,s})_{s\in S},P_4,P_5,$
\item $\qquad \qquad \qquad \qquad \pi^0_0,\dots,\pi^0_{k-1},\pi^1_0,\dots,
\pi^1_k,F_2,F_3,F_4,F_5,(F_{3,g^*})_{g^*\in G} \rangle$
\item and let $\tau = \tau^- \cup \{ Q_s|s\in S\}$.
\end{description}

\end{definition}

Notice that $|\tau| = |G_\l|+|S|+\aleph_0 = 2^\l$,
since $|G_\l| =2^\l$.

\begin{definition}[The sentence
$\psi_k^\lambda$]
The sentence $\psi_k^\lambda \in L_{(2^\l)^+,
    \w}(\tau)$ is the conjunction
\[ \psi_k^\lambda \equiv \bigwedge T_0 \wedge \psi_G \wedge
  \psi_{\bbz_2}\wedge \psi_H \]
  of the first order theory $T_0$ of $M_I$ (for infinite $I$) and the infinitary sentences
  \begin{itemize}
      \item $\psi_G \equiv \forall z_1z_2 ([P_5(z_1)\wedge P_5(z_2)\wedge F_5(z_1)=F_5(z_2)] \to \bigvee_{g^*\in G}F_{3,g^*}(z_1)=z_2)$,
      \item $\psi_{\bbz_2}\equiv \forall y(P_2(y)\leftrightarrow \bigvee_{s\in S}P_{2,s}(y))$,
      \item $\psi_H \equiv \forall y(P_3(y)\leftrightarrow \bigvee_{s\in S}P_{3,s}(y))$.
  \end{itemize}
\end{definition}

We describe in more detail some parts of the previous definition.

\begin{description}
\item $\psi_G$ says that $G$ acts transitively (through the
functions $F_{3,g^*}$) on copies of $G$ (fibers of $P_5$).
\item $\psi_H$ says that there are no ``non-standard fibers'' in $P_2$: every element of $P_2$
    is in some $P_{2,s}$.
\item $\psi_{\bbz_2}$ says that there are no ``non-standard fibers'' in
  $P_3$: every element of $P_3$ is in some $P_{3,s}$
\item Note that, although there are $2^{2^\l}$ sentences in
the logic, we are only using $2^\l$ of them, as witnessed by
$|G|=2^\l$.
\end{description}

We will also use the following variant on the standard model:
for a set $I$ and a function
\[ f: \ialakmasuno\times S \to \bbz_2, \]
we will now build models $M_{I,f}$ and $M^-_{I,f}$.

\begin{definition}\label{defMIf}[The models $M_{I,f}$ and
    $M^-_{I,f}$]
    Let $f: \ialakmasuno\times S \to \bbz_2$ and $I$ a set. Then 
    $M_{I,f}$ is the $\tau$-model constructed just like $M_I$, with only one difference: the
    interpretation of $Q_s$, for $s\in S$, now is  the set  of
tuples  
\[ \langle a_0,\dots,a_k,u_0,\dots,u_k,x_0,\dots,x_{k-1},y_k,z \rangle \]
(see page~\pageref{defM_I}) with condition $(\z)$ replaced by
\[ (\z)_f^* \qquad \bbz_2 \models \sum_{\ell<k} i_\ell =
h_k(u_0) + g(s) + f(u,s). \]
The $\tau^-$-model $M^-_{I,f}$ is then defined as $M_{I,f}\restriction \tau^-$.
\end{definition}

We will use the models $M_{I,f}$ later as canonical ways of describing variants in the choices of
elements of the groups when studying models of the sentence $\psi^\lambda_k$.

We call a $\tau$-structure {\bf strongly standard} if
  $M\restriction \tau^- = M_I\restriction \tau^-$ for $I = P_0^M$.

\begin{definition} {\bf (Some abstract classes related to
    $\psi_k^\lambda$)}
\end{definition}
\begin{enumerate}
\item Let $K_1 := \{ M | M\approx M_{I,f}$ for some infinite set $I$, for
some $f$ as in \ref{defMIf}$\}$. Then $K_1$ is a
class of $\tau$-models.  
\item Let $K^*(\lambda,k) := \Mod(\psi_k^\lambda)$ with the strong
  substructure relation
  \[\prec_{K^*(\lambda,k)}:=\prec_{L_{(2^\lambda)^+,\omega}}.\]
\item $M$ from $K^*(\lambda,k)$ is {\bf standard} if $P_{1,1}^M =
  [P_0^M]^k$ and $P_{1,2}^M = [P_0^M]^{k+1}$ and the $\pi^t_\ell$'s
  correspond to the actual projections sending $u\in \ialak$ to its
  $\ell$'th coordinate in $I$.

\end{enumerate}

\begin{claim}
 For any $M_I\models \psi_k^\lambda$, $M_I \approx M_{I,{\bf 0}}$, for the
 function ${\bf 0}:\ialakmasuno \times S\to \bbz_2$ of constant value
 $0$.
\end{claim}

The proof is immediate from the definition.

\begin{claim}\label{allisostrstrd}
Every $N\models \psi_k^\lambda$ is isomorphic to a strongly standard $M$.
\end{claim}

\Proof Let $N\models \psi_k^\lambda$ and let $I:=P_0^N$. Then $N\restriction \tau^-\approx
M_I\restriction \tau^-$ (following the definition of the sorts of the vocabulary $\tau^-$). Then
\emph{define} the interpretations of the relevant predicates $Q_s$ on $N$ by mapping directly from
their definition on the strongly standard model $M_I$. \qed

\medskip
\noindent
Next, a straightforward observation. 

\begin{claim}
$M_{I,f}$ is strongly standard.
\end{claim}




\begin{proposition}
  $(\cK^*(\lambda,k),\prec_{K^*(\lambda,k)})$ is an abstract elementary
  class with Löwenheim-Skolem number $2^\lambda$.
\end{proposition}

%
%

We do not investigate properties of this AEC in this paper; however, we propose some conjectures at
the end of the paper on their properties and on their connection with good frames and the work of
Boney and Vasey~\cite{BoVa}.

\section{Categoricity of $\psi^\lambda_k$ below $\lambda^{+k}$}

In this section we study the categoricity spectrum of  $\psi^\lambda_k$. The strategy consists of
the following steps:

\begin{itemize}
    \item Since the complexity of models of  $\psi^\lambda_k$ hinges on the predicates $Q_s$,
	and these ultimately depend on \emph{choices} of elements of the copies of the groups above
	the ``supports'' (in the standard case, $k$-element subsets of $(k+1)$-sets of the index set),
	we will develop a language of \textbf{choice functions} to deal with these.
    \item Furthermore, comparing different models will amount to dealing with \textbf{correction
	functions} associated to the choice functions. We also set up a language for these.
    \item Later (Lemma~\ref{labase}) we establish that for every model $N$ of $\psi^\lambda_k$ and global
	choice for $N$ with correction function $f$ there are an index set $I$ and an isomorphism
	${\mathbf h}$ between $N$ and $M_{I,f}$.
    \item Therefore, establishing categoricity in a cardinality $\kappa$ amounts to showing that for
	every $N\models \psi^\lambda_k$ of cardinality $\kappa$ there is a global choice for $N$
	with correction function $0$ (for $\kappa<\lambda^{+k}$; see
	Theorem~\ref{thmchoicecorrectionzero}).
    \item The rest of the section is devoted to showing that if $M\models \psi^\lambda_k$ and
	$|M|<\lambda^{+k}$ then there is a global choice function for $M$ with correction function
	$0$ - with the cardinality restriction in place, we may conclude that $\psi^\lambda_k$ is
	categorical in $\lambda^+,\lambda^{++},\dots,\lambda^{+k-1}$. This part requires several
	lemmas on extending choice functions while keeping the correction function $0$; these
	lemmas depend crucially on the cardinality being of the form $\lambda^{+m}$ for $m$ a
	natural number \emph{below} $k$. This is why the proof in this section only provides
	categoricity up to $\lambda^{+k-1}$.
\end{itemize}

\subsection{Solutions, choices and correction functions}
\label{sec:solch}

We will now define choice functions and correction functions. These will be used
to study models of $\psi^\lambda_k$ of cardinality $\l^+,\dots,\l^{+k-1}$.



Expanding choices from partial to global ones is the crucial issue.

\begin{definition}[Partial $M$-$(J_0,J_1,J_2)$-choice]
    \label{partialchoice}
    For $M\models \psi_k^\lambda$, we say $\chxyz$ is a 
{\bf partial $M$-$(J_0,J_1,J_2)$-choice} if
\begin{description}
\item[(a)] $J_0, J_1\subset P^M_{1,1}$, $J_2\subset P^M_{1,2}$,\\
(so, in the case of standard models, $J_0,J_1\subset
\ialak$, $J_2\subset  \ialakmasuno$)
\item[(b)] $\bar{x} = \langle x_{u,s} | s\in S, u\in J_0\rangle$,
  where
  \[ x_{u,s}\in (P_{3,s}^M)^{-1}(u)\subset P_{3,s}^M.\]
\item[(c)] $\bar{y} = \langle y_{u,s} | s\in S, u\in J_1\rangle$,
\[ y_{u,s}\in H^M_{u,s} :=  (P_{2,s}^M)^{-1}(u) \subset P_{2,s}^M.\]
\item[(d)] $\bar{z} = \langle z_u | u\in J_2\rangle$,
\[ z_u\in G_u^M := (F_5^M)^{-1}(u) \subset P_5^M.\]
\end{description}
\end{definition}

Therefore $\bar{x}$ essentially chooses an element $i$ in the corresponding copy of
$\bbz_2$, $\bar{y}$ chooses a $h$ in the corresponding copy of $H_I$, $\bar{z}$
chooses a $g$ in the corresponding copy of $G$, for each relevant $(u,s)$.

So,
$x_{u,s}$ is {\it some} element in the `fiber' of $u$ via $F_3^M$, and
analogously for $\bar{y}$ and $\bar{z}$.

\begin{definition}
    We call  $\chxyz$ a
{\bf partial $M$-$J$-choice} if it is an $M$-$(J,J,J^M_*)$-choice, where

\[J_*^M :=
\Big\{ a \in P_{1,2}^M \Big| \bigwedge_{m\leq k} \exists b\in J
       [\bigwedge_{\ell<m} 
(\pi^1_\ell(a) = \pi^0_\ell(b) \wedge \bigwedge _{\ell\in [m,k[}
      \pi^0_\ell(b) =
\pi^1_{\ell+1} (a)]\Big\} .\] 
 Similarly, we say that $\chxyz$ is a {\bf global $M$-choice} if it is a
partial $M$-$P_{1,1}^M$-choice. We will sometimes just say
``$M$-choice'' (if clear from context).
   \label{partialMJchoice}
\end{definition}

\noindent
The previous is a way of describing, in our language of projections,
that (in the standard case) $J_*^M$ consists of the $k+1$-element sets
such that \emph{all} their ($k+1$-many) $k$-element subsets are in
$J$).

So, when $M$ is standard, we have that
\[J^M_* = \Big\{ \langle a_\ell | \ell \leq k \rangle \Big| \bigwedge
_{m\leq k} \langle a_\ell |\ell \not= m \rangle \in J \Big\} .\]


\begin{definition}
    Fix a standard $M$ and a
$M$-$(J_0,J_1,J_2)$-choice $\chxyz$. Then we let the {\bf correction
function} $f$ for $M$ and $\chxyz$ be the function such that
\begin{enumerate} 
\item $\Dom (f)$ is the set of pairs
$(u,s)$ such that
\begin{description}
\item[$(\a)$] $u = \langle a_\ell | \ell\leq k\rangle \in J_2\subset
P^M_{1,2}$,
\item[$(\b)$] if $u_m := \langle a_\ell | \ell\leq k, \ell\not=
  m\rangle$,
$u_\ell \in J_0$ for $\ell<k$, $u_k\in
J_1\subset P^M_{1,1}$,
\end{description}
\item $\ran(f) \subset \bbz_2$, and 
\item (recall $x_{u_\ell,s},y_{u_k,s},z_{u_k}$ are from the choice)
\[ f(u,s) = 0 \Leftrightarrow \langle
a_0,\dots,a_k,u_0,\dots,u_k,x_{u_0,s},
\dots,x_{u_{k-1},s},y_{u_k,s},z_{u_k}\rangle
\in Q_s^M.\]
\end{enumerate}
    \label{correctionfunction}
\end{definition}


The next claim is a general observation on correction functions
and choices.


\begin{claim}
For every $M\in \Mod(\psi_k^\lambda)$, there is an $M$-choice
$\chxyz$.
\end{claim}

\Proof Immediate: just construct the tuples. There the demands are on
each choice separately. There are no demands connecting different
choices. \qed

The next lemma is a crucial step. It shows how to build possible
isomorphisms from an arbitrary model $N$ of $\psi^\lambda_k$ to standard models
$M_{I,f}$.

\begin{lemma}\label{labase}
    Let $N\in \Mod(\psi_k^\lambda)$ and let $\chxyz$ be a global $N$-choice
    with correction function $f$. \underline{Then,} there exist a set $I$ and an isomorphism
    \[ {\mathbf h}:N\to M_{I,f}.\]
    Furthermore, the isomorphism behaves as follows on the global $N$-choice $\chxyz$:
\[ {\mathbf h}(x_{u,s}) = ({\mathbf h}(u),s,0_{\bbz_2}), \quad {\mathbf h}(y_{u,s}) =
({\mathbf h}(u),s,0_{H_I}), \quad {\mathbf h}(z_u) = ({\mathbf h}(u),0_G) .\]
\end{lemma}

\Proof Let $N\models \psi_k^\lambda$, and fix a global $N$-choice
$\chxyz$ with correction function $f$. We build $I$ and $\mathbf
h$ as in the statement.

First, we extract the predicates for the model $M=M_{I,f}$: let
$I:=P^N_0$. Clearly,
$P^M_0 = P^N_0$.

We now define $\mathbf h$, following the predicates of the domain of $N$ (remember that the domain
of $N$ is the disjoint union $$P_0^N\cup P_{1,1}^N\cup P_{1,2}^N\cup P_2^N\cup P_3^N\cup P_4^N\cup
P_5^N$$ and the predicates $P_2^N$ and $P_3^N$ are each partitioned into classes $P_{2,s}^N$,
$P_{3,s}^N$ ($s\in S$)).
\begin{itemize}
    \item $\mathbf h$ is the identity on $P^N_0=I=P^M_0$.
    \item if $x\in P_{1,1}^N$, $\ell<k$, $\pi^0_\ell(x)=x_\ell$ ($\in P_0^N$), then ${\mathbf
	h}(x):=({\mathbf h}(x_0),\dots,{\mathbf h}(x_{k-1}))$.
    \item similarly, if $x\in P_{1,2}^N$, $\ell<k+1$, $\pi^1_\ell(x)=x_\ell$ ($\in P_0^N$), then
	${\mathbf h}(x):=({\mathbf h}(x_0),\dots,{\mathbf h}(x_k))$.
    \item if $x\in P_{2,s}^N$ then ${\mathbf h}(x)=({\mathbf h}(F_2^N(x)),s,-)\in [I]^k\times S\times
	H_I$. For now we only know the third coordinate must be an element of $H_I$.
	Also, \textbf{as soon as we know} the third coordinate of the image of one
	element $x_0$ of a fiber inside the predicate $P_{2,s}^N$, we also know the third coordinate for all other
	elements $x$ of that fiber: since the action given by $F_4^N$ is transitive (as encoded
	by $T_0$), there is some $h_0\in P_4^N$ such that $F_4^N(x_0,h_0)=x$. Then (if we also have
	a definition of $\mathbf h$ on elements of $P_4^N$), we have that ${\mathbf h}(x)={\mathbf
	h}(F_4^N(x_0,h_0))=F_4^N({\mathbf h}(x_0),{\mathbf h}(h_0))$.
    \item Similarly, if $x\in P_{3,s}^N$, then ${\mathbf h}(x)=({\mathbf h}(F_3^N(x)),s,-)\in
	[I]^k\times S\times \bbz_2$ and just as before the value of $\mathbf h$ on one element of
	the fiber will determine the rest.
    \item And similarly, if $x\in P_5^N$, then ${\mathbf h}(x)=(F^N_5(x),-)\in [I]^{k+1}\times G$.
    Again, since $N\models \psi_G$, the action (``of $G$'') encoded by the family of functions
    $F_{3,g^*}^N$ is transitive, and therefore knowing a second coordinate for \emph{one} element of
    a fiber of 
    $P_5^N$ implies knowing it for all elements of the corresponding fiber that predicate.
\end{itemize}

It therefore remains, in order to complete the definition, to make \textbf{choices} of images of
elements of $P_4^N$ (images in $H_I$ - but this is easy, as $H_I$ is definable in our structure) and selecting, for each $s\in S$, one image in each one of the
relevant fibers. We now use the correction function $f$ and the predicates $Q_s$.

So fix $s\in S$. Checking the equivalence we are looking for, namely
\[ Q^N_s(a_0\dots a_ku_0\dots u_kx_0\dots x_{k-1}y_kz)\]
\[ \Updownarrow\]
\[ Q^{M_{I,f}}_s({\bf h}(a_0)\dots {\bf h}(a_k){\bf h}(u_0)\dots {\bf
h}(u_k){\bf h}(x_0)\dots {\bf h}(x_{k-1}){\bf h}(y_k){\bf h}(z)),\]

amounts to answering the question

\[ \bbz_2 \models \sum_{\ell<k}i_\ell = h_k(u_0) + g(s) + f(u,s) \]
\[ \Updownarrow ?\]
\[ \bbz_2 \models \sum_{\ell<k}{\bf h}(i_\ell) = h_k({\bf h}(u_0)) +
g({\bf h}(s)) + f({\bf h}(u),{\bf h}(s)) \]

\noindent
Now letting 

\medskip

$\qquad \left \{ \begin{tabular}{l}
${\bf h}(x_{u,s}) = ({\bf h}(u),s,0_{\bbz_2})$,\\
${\bf h}(y_{u,s}) =
({\bf h}(u),s,0_H),$\\
${\bf h}(z_u) = ({\bf h}(u),0_G)$
\end{tabular} \right .$

\medskip
\noindent
works for these equations: we are assigning $0$ on the missing
coordinates (third or second) -- exactly to those elements of the fibers ($x_{u,s}$, $y_{u,s}$,
$z(u)$) that had already been picked by the choice function.

\noindent
Why is this enough?

\noindent
Well, our definition turns the equation (at the choices) into
\[ \bbz_2 \models 0=\sum_{\ell<k}0 = 0(\star) + 0(\star) +f(\star).\]
But, since $f$ was a correction function \emph{for the choice function} $\chxyz$,
\[ f(u,s) = 0 \Leftrightarrow \langle
a_0,\dots,a_k,u_0,\dots,u_k,x_{u_0,s},
\dots,x_{u_{k-1},s},y_{u_k,s},z_{u_k}\rangle
\in Q_s^N,\]
and therefore our definition of $\bf h$ works.
\qed


\begin{definition}[Canonical choice]
Fix $M=M_{I,f}$, and let $\chxyz$ be the $M$-choice
given by
\[ x_{u,s} = (u,s,0_{\bbz_2}),\]
\[ y_{u,s} = (u,s,0_{H_I}),\]
\[ z_u = (u,0_G).\]
This is by definition the \textbf{canonical $M$-choice}.
\end{definition}

\begin{claim}\label{zerosforchoices}
\begin{enumerate}
\item If $\chxyz$ is a global $M$-choice, $M\models \psi_k^\lambda$, and $f$ is
the $M$-correction function for $\chxyz$, and $f$ is identically zero,
then $M\approx M_I$ for some $I$.
\item If $f$ above is zero on $P_{1,1}^M$, $P_{1,2}^M$ and
$f=f'\restriction J_2\times S$, then $M\approx M_{P_1,f'}$.
\end{enumerate}
\end{claim}

\Proof Part (1) is a consequence of \ref{labase}. Part (2) is
clear.
\qed


\begin{corollary}
The correction function for $M_{I,f}$ with the canonical $M$-choice
$\chxyz$ is $f$.
\end{corollary}

\Proof Similar to the previous: add zeroes to $f$ as
in~\ref{zerosforchoices}. \qed

\subsection{Models of cardinality below $\lambda^{+k}$}

The rest of the section contains several \emph{extension lemmas} for models of $\psi^\lambda_k$ of
cardinalities $\lambda,\lambda^+$, etc.: the crucial issue is to build a choice function with
null correction function. This may be started first at cardinality $\lambda$, and then pushed up.
But each step up exacts an ``amalgam of choices'' possible only up to cardinality
$\lambda^{+k-1}$.

The next lemma is the first step in the categoricity proof. It provides a specific kind of
extension of choice: from an $M$-$J$-choice with correction function zero to a global $M$-choice
with correction function zero \emph{when $J$ consists of $k$-subsets of the ``support part'' of $M$
that \emph{omit} some fixed set $W$ of \textbf{at most} $k$-many elements.} Also, it is worth
stressing the lemma is about standard $M$.

For instance, when $m=2<k=3$, the lemma would mean we start with a choice function $\chxyz$ for all
``triangles'' and ``tetrahedra'' \emph{omitting} some fixed pair $\left\{ a,b \right\}$ \dots and
then would extend the choice function (with correction function zero) to all triangles and
tetrahedra.

\begin{lemma}\label{existxyz}[Extension property for $W$ of size
    $m<k$, $|P_0^M|\leq \l$]\\
Assume $m<k$, $M\models \psi^\lambda_k$, $M$ is strongly standard, $|P_0^M|\leq
\l$, 
$W\subset P_0^M$, $W=\{ b_\ell |\ell < m\}$ with no
repetition, $J=\{ u\in P^M_{1,1} |  W\not\subset u\}$
(note that $u\in [P^M_0]^k$, as
$M$ is standard), $\chxyz$ is an $M$-$J$-choice
with correction function $f_0$, identically zero.
{\sf Then,} we can extend $\chxyz$ to an $M$-choice
with correction function identically zero.
\end{lemma}

\Proof
\begin{description}
\item[Part A:]
Without loss of generality, by~\ref{allisostrstrd}, since $M$ is
strongly standard, $I=P^M_0$.
Let $\langle \bar{a}^\a | \a<\b^*\rangle$ list $P_{1,1}^M$ with $\langle
\bar{a}^\a | \a < \a^*\rangle$ listing $J$ (we have also used $u$ for
naming these $\bar{a}^\a$'s).
Let $\langle \bar{b}^\g | \g<\g^*\rangle$ list $\{ \bar{a}\in \ialakmasuno |
\bar{a}$ with no repetition and $W\subset \ran (\bar{a})\}$ and $\g^*
<\l^+$.

Our hypothesis is then that we have choice functions for all $u\in
P_{1,1}^M$ such that
$u\not\supset W$, with correction function zero.

We list these choice functions as follows:
Let, for $\a <\a^*$,
\begin{description}
 \item $x_{\bar{a}^\a,s} = (\bar{a}^\a,s,i_{\a,s})\in
 (\bbz_2)_{\bar{a}^\a,s}$, $i_{\a,s}\in \bbz_2$,
\item $y_{\bar{a}^\a,s} = (\bar{a}^\a,s,h_{\a,s})\in
 H_{\bar{a}^\a,s}$, $h_{\a,s}\in H_I$,
\item $z_{\bar{b}^\g} = (\bar{b}^\g,g^\g)$, $g^\g \in G$.
\end{description}

We now have to extend these choice functions to those $u$ such that $u\supset W$.

We will now choose $x_{\bar{a}^\a,s} = (\bar{a}^\a,s,i_{\a,s})$,
$y_{\bar{a}^\a,s} = (\bar{a}^\a,s,h_{\a,s})$, $z_{\bar{b}^\g} =
(\bar{b}^\g,g^\g)$ for $\a^* \leq \a <\b^*$ and appropriate $\g$.

Without loss of generality, $\b^*\leq \a^*+\l$, $\g^*\leq
\l$. (Remember $S=[\l]^{<\aleph_0}$.)

\item[Part B:]
\begin{description}
\item[First,] we choose $i_{\a,s}=0_{\bbz_2}$ for $\a^* \leq \a < \b^*$,
    $s\in S$. This provides the choices $x_{\bar{a}^\alpha,s}$ for $\alpha^*\leq \alpha<\beta^*$.
\item[Second,] we choose the relevant $h$ functions. We try a value for $h_{\a,s}$ for $\a^* \leq \a < \b^*$
and $s\in S$ so that
\begin{description}
\item[(*)] if $\g \in s\subset \l$, $\bar{b}^\g = \langle b^\g_\ell |
\ell\leq k \rangle$, 
$u_n^\g = \langle b^\g_\ell | \ell \leq k, \ell\not= n\rangle$, let
$\varepsilon(\gamma,n) < \b^*$ be such that $u^\gamma_n =
\bar{a}^{\varepsilon(\g,n)}$ then
\end{description}
\[ h_{\varepsilon(\g,k),s}(\bar{a}^{\varepsilon(\g,0)}) = 0\]
\[\Updownarrow \]
\[ \langle b_0^\g,\dots ,b_k^\g,u_0^\g,\dots,u_k^\g,x_{\varepsilon(\g,0),s},
    \dots,x_{\varepsilon(\g,k-1),s},(\bar{a}^{\varepsilon(\g,k)},s,0_H),
(u^\g,0_G)\rangle \in Q^M_s .\]

\end{description}
Note that all the elements in the bottom part of the previous have already been defined previously. 

\noindent
Let $t(\gamma,s)$ be $0$ if the bottom statement is true, $1$
otherwise.
For our fixed $s\in S$, let $A_s$ be the (finite) set $\left\{ \varepsilon(\gamma,k)\mid \gamma\in s
\right\}$; we now define $h_{\alpha,s}$ for our fixed $s$ and at the relevant $u$. If
$\alpha\notin A_s$, then let $h_{\alpha,s}(u)=0$ for all $u$. If $\alpha\in A_s$, we proceed as
follows. First we consider the set
$s_\alpha:=\left\{ \gamma\in s\mid \varepsilon(\gamma,k)=\alpha \right\}$ and we then define
\[ h_{\alpha,s}(u)=\begin{cases}
	t(\gamma,s), & \mbox{ if }u=a^{\varepsilon(\gamma,0)}\mbox{ for some }\gamma\in s_\alpha,\\
	    0, & \mbox{ otherwise.}
\end{cases} \]


Notice that these decisions are made \textbf{for each $s$ separately,} and that as we fix $s$ we
really deal with one $\alpha\in [\alpha^*,\beta)$: when we choose $h_{\alpha,s}$ we only have to
consider $\gamma<\gamma^*$ such that $\varepsilon(\gamma,\ell)\in s$. There are only finitely many
such $\gamma$'s. Moreover, if $\gamma_1\not= \gamma_2\in s$ and
$\varepsilon(\gamma_1,k)=\alpha=\varepsilon(\gamma_2,k)$ then necessarily
$\varepsilon(\gamma_1,0)\not=\varepsilon(\gamma_2,0)$, as $\bar{b}^\gamma$ is reconstructible from
$\alpha$ and $\varepsilon(\gamma_1,0)$. 

So, our definition of the functions $h_{\varepsilon(\gamma,k),s}$ does not have contradictory
demands; since the set $s_\alpha$ is finite, the function defined has finite support.



\item[Part C:]
    Having extended the choices $x$ and $y$, it only remains to extend the $z$ part.
    Let us now fix $\gamma$ and find a $g\in G$ that will provide a choice (with correction function
    zero) for the corresponding $\bar{b}^\gamma$. [Recall that if
    $\bar{b}\in \ialakmasuno$ is such that $\bar{b}\supset W$, then
$\bar{b}=\bar{b}^\g$ for some $\g<\g^*$.]

    But then the set
\[ S^*_\g = \Big\{ s\in S \Big| M\models Q_s\big( b_0^\g,\dots,b_k^\g,
	\bar{a}^{\epsilon(\g,0)},\dots,\bar{a}^{\epsilon(\g,k)},\qquad \qquad \qquad \qquad \]
	\[ \qquad \qquad \qquad \qquad \qquad
x_{u^\g _0,s},\dots,x_{u^\g _{k-1},s},
y_{u^\g _k,s},(u^\g,0_G)\big) \Big\} \]
belongs to $\frak D$. This last point holds by the regularity of $\frak D$:
if $s_0\in S^*_\gamma$ then
the tuple $$\big( b_0^\g,\dots,b_k^\g,
	\bar{a}^{\epsilon(\g,0)},\dots,\bar{a}^{\epsilon(\g,k)},x_{u^\g _0,s_0},\dots,x_{u^\g
	_{k-1},s_0},
y_{u^\g _k,s_0},(u^\g,0_G)\big)$$ belongs to $Q_{s_0}$; now, if $s\supset s_0$, the corresponding tuple
$$\big( b_0^\g,\dots,b_k^\g,
	\bar{a}^{\epsilon(\g,0)},\dots,\bar{a}^{\epsilon(\g,k)},x_{u^\g _0,s},\dots,x_{u^\g _{k-1},s},
y_{u^\g _k,s},(u^\g,0_G)\big)$$ will belong to $Q_s$.

Next choose $z_{\bar{b}^\gamma} := (\bar{b}^\gamma,g)$ with $g$ given by

\medskip
\qquad $g(s) = \left \{ \begin{tabular}{ll}
$0$	&	if $s\in S^*_\g$\\
$1$	&	if $s\not\in S^*_\g$
\end{tabular} \right . $

\medskip
Now then, with these $x$, $y$ and $z$, the equation holds.
\end{description}
\qed$_{\ref{existxyz}}$

We now deal with \textbf{systems} of choices, trying to obtain extensions with correction function
zero at cardinalities above $\lambda$. In what follows, as usual, ${\mathcal P}^-(m_2)$ denotes
${\mathcal P}(m_2)\setminus \left\{ m_2 \right\}$. 

\begin{definition}[Compatible system of choices]
    Let $M\models \psi^\lambda_k$ be strongly standard, $A_\emptyset \subset
P_0^M$, $m_1+m_2<k$ and
$a_0, \dots,a_{m_2-1}$ different elements of $P_0^M\setminus
A_\emptyset$. Then
$$\langle A_s,\chxyz_s | s\in {\cal P}^-(m_2) \rangle$$
is a \textbf{compatible
$\l^{+m_1}$-${\cal P}^-(m_2)$-system of choices} iff
\begin{enumerate}
\item $\bigcup_{s\in {\cal P}^-(m_2)}A_s = A_\emptyset \cup \{
a_0,\dots,a_{m_2-1} \}$, $|A_\emptyset| \leq \l^{+m_1}$, $A_s
= A_\emptyset \cup \{ a_t | t\in s\}$.
\item $\chxyz_s$ is a $M$-$[A_s]^k$-choice, for each $s\in {\cal P}^-(m_2)$.
\item For every $s,t\in {\cal P}^-(m_2)$, $s\subset t \Rightarrow \chxyz_s
\subset \chxyz_t$\footnote{here, of course, we are abusing notation - by
$\chxyz_s \subset \chxyz_t$ we mean $\bar{x}_s\subset \bar{x}_t$,
$\bar{y}_s\subset \bar{y}_t$ and $\bar{z}_s\subset \bar{z}_t$.}.
\end{enumerate}
\end{definition}

\begin{lemma}\label{0case}
If $\langle A_s,\chxyz_s | s\in {\cal P}^-(m_2)\rangle$ is a
compatible $\l$-${\cal P}^-(m_2)$-system with $m_2<k$ (with correction function zero for each $s\in
{\cal P}^-(m_2)$),
\underline{then} there is an
$M$-$\bigcup_{s\in {\cal P}^-(m_2)} A_s$-choice $\chxyz$ extending all
the $\chxyz_s$, for $s\in \pmin(m_2)$, with correction function zero.
\end{lemma}

\Proof Let $m_2<k$ and let $\langle A_s,\chxyz_s | s\in {\cal P}^-(m_2)\rangle$ be a compatible
$\l$-${\cal P}^-(m_2)$-system, each choice in the system with correction function zero.
Notice that
\[ u\in [\bigcup_{s\in {\cal P}^-(m_2)} A_s]^k \setminus
\bigcup_{s\in {\cal P}^-(m_2)}[A_s]^k,\]
if and only if $\{ a_0\dots a_{m_2-1}\} \subset u$.

\noindent
We first notice that by compatibility, the union of the choices $\chxyz_s$ along ${\mathcal
P}^-(m_2)$ is an $M$-choice for $\bigcup_{s\in {\cal P}^-(m_2)}[A_s]^k$. It remains to extend that
choice to an $M$-$\bigcup_{s\in {\cal P}^-(m_2)} A_s$-choice $\chxyz$ with
correction function zero.

\noindent
We may apply Lemma~\ref{existxyz} (here, the set $W$ of cardinality $m_2<k$ is
    $\left\{ a_0,\dots,a_{m_2-1} \right\}$ and
the lemma provides the extension from a $M$-$\bigcup_{s\in {\cal P}^-(m_2)}[A_s]^k$-choice with
correction function zero to a $M$-$\bigcup_{s\in {\cal P}^-(m_2)} A_s$-choice $\chxyz$ with
correction function zero - we extend the choice from those $k$-sets omitting $W$ to all of them).
\qed$_{\ref{0case}}$

%

\begin{lemma}\label{gencase}
    Let $m_1+m_2<k$. 
If $\langle A_s,\chxyz_s | s\in \pmin(m_2)\rangle$ is a compatible
$\l^{+m_1}$-$\pmin(m_2)$-system of choices with correction function zero, \underline{then} there is
a $\bigcup_{s\in \pmin(m_2)}A_s$-choice $\chxyz$ with correction function zero such that $\chxyz_s
\subset \chxyz$ for every $s\in \pmin(m_2)$.
\end{lemma}

\Proof By induction on $m_1$. For $m_1=0$, this is lemma \ref{0case}. For
$m_1 >0$, suppose $A_s = A_\emptyset \cup \{ b_j|j\in s\}$. Enumerate
$A_\emptyset$ as $\langle a_\b | \b <\l^{+m_1}\rangle$. Let
$A_\emptyset^\a = \{ a_\b |\b<\a\}$ and $A_s^\a = A_\emptyset^\a
\cup \{ b_j|j\in s\}$ for every $s\in \pmin(m_2)$. Finally, let $\chxyz_s^\a$ be
the restriction of the choice we have $\chxyz_s$ from the compatible system (with correction
function zero) to an $M$-$A^\a_s$-choice (also immediately with correction function zero).

The plan is to obtain an $M$-$\bigcup_{s\in m_2}A^\alpha_s$-choice $\chxyz_\a$ with correction function zero
for each $\alpha<\lambda^{+m_1}$, such that $\a <\beta < \l^{+m_1}$ implies
$\chxyz_\alpha\subset \chxyz_\beta$.

We build $\chxyz_\alpha$ by another induction, on $\alpha<\lambda^{+m_1}$. For $\alpha=0$, the empty
choice function is an $M$-$\emptyset$-choice ($A^0_s=\emptyset$ for each $s$). When $\alpha$ is a
limit ordinal, the union of the chain of choices $(\chxyz_\beta)_{\beta <\alpha}$ is a $M$-$\bigcup_{s\in
m_2}A^\alpha_s$-choice with correction function zero. Finally, for $\alpha=\beta+1$, we proceed as
follows: we already have, by induction hypothesis, an $M$-$\bigcup_{s\in {\mathcal
P}^-(m_2)}A^\beta_s$-choice with correction function zero, $\chxyz_\beta$; consider also the choices
$\chxyz_s^\alpha$ for $s\in {\mathcal P}^-(m_2)$. Since the cardinalities of all their domains are
$<\lambda^{+m_1}$, we may without loss of generality regard the previous choices as forming a
compatible $\lambda^{+m_1-1}$-${\mathcal P}^-(m_2+1)$-system of choices with correction function
zero: the set $\left\{ b_i\mid i\in s\right\}\cup \left\{ \beta \right\}$ has cardinality $m_2+1$.
Since $(m_1-1)+(m_2+1)=m_1+m_2<k$, we may apply the induction hypothesis; we obtain $\chxyz_\alpha$
an $M$-$\bigcup_{s\in {\mathcal P}^-(m_2)}A_s^\alpha$-choice with correction function zero.

Having constructed this chain $\chxyz_\alpha$ for $\alpha<\lambda^{+m_1}$, we just let
\[ \chxyz:=\bigcup_{\alpha<\lambda^{+m_1}}\chxyz_\alpha.\]
This is a $\bigcup_{s\in {\mathcal P}^-(m_2)}$-choice with correction function zero, extending all
the choices in the system.

\qed$_{\ref{gencase}}$

We may now obtain our general extension property.

\begin{lemma}\label{fullext2002}{\bf (Full extension)}\\
Let $M\models \psi$ be strongly canonical,
$J_1\subset J_2\subset P_0^M$, with $|J_2| <
\l^{+k-1}$ and $\chxyz$ an $M$-$J_1$-choice with correction function
identically zero. Then $\chxyz$ can be extended to an $M$-$J_2$-choice
with correction function identically zero.
\end{lemma}

\Proof Without loss of generality, $J_2 = J_1 \cup \{ b\}$. If $J_1$
has size $\leq \l$,
this
is lemma \ref{existxyz}. Now suppose $|J_1| = \l^{+m_1}<\lambda^{+k-1}$ (therefore
$m_1<k-1$) and enumerate $J_1$ 
as $\langle a_\b | \b<\l^{+m_1}\rangle$. Let $J_1^\a = \{
a_\b | \b<\a\}$, and let $\chxyz_\a$ be the restriction
of $\chxyz$ to an $M$-$J_1^\a$-choice. We define by induction
$M$-$J^\a_1$ choices with correction function identically zero
$\chxyz_\a' \supset \chxyz_\a$.

We may use here lemma \ref{gencase} for $m_2 = 2$
to extend $\chxyz_\a' \cup \chxyz_{\a+1}$ to an $M$-$J_1^{\a+1}\cup\{
b\}$-choice with correction function identically zero: since $m_1<k-1$, along the induction the
cardinality is $<\lambda^{+m_1}$, say $\lambda^{+m_1'}$ for some $m_1'<m_1$. Since we also have that
$m_1<k-1$, then $m_1'+2<k$ and we can use $m_2=2$ when invoking lemma~\ref{gencase}.

At limits take
unions; finally,
\[ \Big(
\bigcup_{\a<\l^{+m_1}}\bar{x}_\a',\bigcup_{\a<\l^{+m_1}}\bar{y}_\a',
\bigcup_{\a<\l^{+m_1}}\bar{z}_\a'\Big) \]
is an $M$-$J_2$-solution extending $\chxyz$.
\qed

\begin{theorem}
If $M\models \psi_k^\lambda$ is strongly standard
and $|M|<\l^{+k}$ then there is an $M$-choice with
correction function identically zero.
\label{thmchoicecorrectionzero}
\end{theorem}

\Proof We apply
Lemma~\ref{fullext2002} (starting from the empty choice function, and taking unions at limits): the
lemma gives an extension of a choice function with correction function zero
from $J_1\subset P_{1,1}^M$ to $J_2$ with $J_1\subset
J_2\subset P_{1,1}^M$ provided $|J_2|<\lambda^{+k-1}$. Here $|M|$ may be \emph{equal to}
$\lambda^{+k-1}$ (at ``worst''); if (in that case) we enumerate $P_{1,1}^M$ as $\left\{ a_\beta\mid
\beta<\lambda^{+k-1} \right\}$ then given $\alpha < \lambda^{+k-1}$, $|\left\{ a_\beta\mid
\beta<\alpha \right\}|<\lambda^{+k-1}$ and we can apply Lemma~\ref{fullext2002} to get an extension
of the choice with correction function zero to $\left\{ a_\beta\mid \beta<\alpha \right\} $.
\qed

\begin{theorem}\label{psicateginlow}{\bf (Categoricity and
amalgamation up to $\lambda^{+k}$)}\\
\begin{enumerate}
\item For $m< k$, $\Mod(\psi_k^\lambda)$ has a unique strongly standard model $M$, $|P_0^M| =
    \l^{+m}$, modulo isomorphism.
\item For $m<k-1$, if $2^\l \leq \l^{+m}$, then $\cK^*(\lambda,k)$ has
amalgamation in $\l^{+m}$.
\item If $m< k$, $\l^{+m} \geq 2^\l$, then $\cK^*(\lambda,k)$ is
  categorical in $\l^{+m}$.
\end{enumerate}
\end{theorem}

\Proof
\begin{enumerate}
    \item Let $N\models \psi^\lambda_k$ be a strongly standard model with $\lambda\leq |P_0^N|<
	\lambda^{+k}$. By Lemma~\ref{labase}, once we have $\chxyz$ a global $N$-choice with
	correction function $f$, then $N\approx M_{I,f}$ for $I=P_0^M$. Now, since $N$ is standard
	and $|P_0^N|\in [\lambda,\lambda^{+k-1}]$, Theorem~\ref{thmchoicecorrectionzero} gives a
	global $N$-choice $\chxyz$ with correction function identically zero (as $N$ is strongly
	standard). So, $N\approx M_I$.
    \item In the proofs of the previous lemmas, amalgamation of choices (along systems) with
	correction function zero is carried out in detail. These give rise to the corresponding
	embeddings and amalgams of models, if the size of these is controlled by the size of their
	$P_0^M$ parts. The only part of a standard model where this cardinality may increase is
	given by the coding of the action of $G$ (remember $|G|=2^\lambda$). If $2^\lambda\leq
	\lambda^{+m}$, the model $M$ will have the same size as $P_0^M$.
    \item Let  $m< k$, $\l^{+m} \geq 2^\l$ and let $M$ be a model in $\cK^*(\lambda,k)$ of size
	$\lambda^{+k}$. Then by Lemma~\ref{allisostrstrd}, $M$ is isomorphic to a strongly standard
	model $N$; also, since $2^\lambda\leq \lambda^{+m}$, $|P_0^M|=\lambda^{+m}$. Thus by part (1)
	$M \approx N\approx M_I$ for some $I$ of size $\lambda^{+m}$.
\end{enumerate}

\qed

\section{Failure of categoricity of $\psi^\lambda_k$ at $\beth_{k+1}(\lambda)$}

We have proved in \ref{psicateginlow} that $\psi$ is categorical in
$\l^{+m}$  if $m<k$ and $2^\l <\l^{+m}$. We now prove that our
sentence is not categorical in any cardinality $\kappa\geq
\mu=\beth_{k+1}(\lambda)^+$. (It is also possible to show that $\psi_k^\lambda$ has the maximal
number of models possible in $\mu$ for each $\mu\geq \beth_{k+1}(\lambda)^+$. We do not do that in
this paper.)

As before, we use our terminology of ``solutions and corrections
functions'' to count the number of models.

\subsection{Combinatorial criteria for (failure of) isomorphism}

In this section we prove a combinatorial criterion for \textbf{non-isomorphism} between two models
of the form $M_{I,f}$.

Before giving the purely combinatorial criterion, we prove the following lemma (a criterion for
isomorphism in terms of \textbf{choices} and \textbf{correction functions}).

\begin{lemma}\label{modset}
If $M_1$ and $M_2$ are strongly standard, and
$\chxyz_\ell$ is an $M_\ell$-choice for $M_\ell$ ($\ell =1,2$),
$P_0^{M_1} = P_0^{M_2}$ with
correction function $f_\ell$ for $\ell = 1,2$ then the following are
equivalent:
\begin{description}
\item[(a)] there is an isomorphism from $M_1$ onto $M_2$ over the
identity on $P_0^{M_1}\cup
P_1^{M_1}$ 
\item[(b)$_1$] there is an $M_2$-choice $\chxyz$
whose correction function is $f_1$,
\item[(b)$_2$] there is an $M_1$-choice $\chxyz$
whose correction function is $f_2$,
\item[(c)] there are functions $g_1$, $g_2$, $g_3$ (``to
correct the choice of zeros''), with
\begin{enumerate}
\item $g_1:\ialak \times S \to \bbz_2$ (like the $x_{u,s}$'s above),
\item $g_2:\ialak \times S \to H_I$ (like the $y_{u,s}$'s above),
\item $g_3:\ialakmasuno  \to G$ (like the $z_u$'s above),
\item if $\langle a_0\dots a_ku_0\dots u_kx_0\dots
x_{k-1}y_k,z\rangle$ are like in Definition \ref{defM_I} for $M_1$,
or $M_2$ \underline{then} 
\[ \bbz_2 \models \sum_{\ell<k} i_\ell - h_k(u_0) - g(s) = \sum
_{\ell<k} g_1(u_\ell,s) - g_2(u_k,s)(u_0) - g_3(u)(s) \]
\end{enumerate}
\end{description}
\end{lemma}

\Proof \begin{description}
\item[(a) $\to$ (b)$_1$] Recall that $M_1\restriction \tau^- =
M_2\restriction \tau^-$, so $M_1$ and $M_2$ have the same universes.
Fix $F:M_1 \stackrel{\approx}{\longrightarrow}_{P_0^{M_1}
\cup P_1^{M_1}} M_2$. 
We have, since $f_1$ is a correction function for $M$ for
the choice $\chxyz_1$, that
\[ f_1(u,s) = 0 \Leftrightarrow
\langle a_0\dots a_k u_0\dots u_k x^1_{u_0s} \dots x^1_{u_{k-1}s} y^1_
{u_ks} z^1_u \rangle \in Q_s^{M_1}.\]
But the right hand side holds iff
\[ \langle a_0\dots a_k u_0\dots u_k F(x^1_{u_0s}) \dots
F(x^1_{u_{k-1}s})
F(y^1_{u_ks}) F(z^1u) \rangle \in Q_s^{M_2},\]
since $F$ is an isomorphism fixing $P_0^{M_1}\cup P_1^{M_1}$, and
$a_0,\dots ,a_k\in P_0^{M_1}$. This gives 
us the $M_2$-choice for which $f_1$ is a correction function: given
$u_\ell \subset u$, $u_\ell\in \ialak$, $u\in \ialakmasuno$, let
$x_{u_\ell,s}' = F(x_{u_\ell,s}^1)$, $y_{u_,s}' = F(y_{u_k,s}^1)$,
$z_u' = F(z_u^1)$.
\item[(a) $\to$ (b)$_2$] Same.
\item[(b)$_\ell$ $\to$ (c)] ($\ell=1,2$) The point of (c) is that we
may find concrete 
representations $g_1,g_2,g_3$, that
act {\it independently from $M_1$ or $M_2$}  as `corrected choice functions'
for the zeros for $f_1$ and $f_2$. So, suppose we have 
a $M_2$-choice $\chxyz$ with correction function $f_1$. Then for any
$u \in P_0^{M_2}$ and any $s\in S$, if $\langle a_0\dots a_ku_0\dots
u_kx_0\dots 
x_{k-1}y_k,z\rangle$ are like in Definition \ref{defM_I} 
\[ \langle a_0\dots a_k u_0\dots u_k x_{u_0s} \dots x_{u_{k-1}s} y_
{u_ks} z_u \rangle \in Q_s^{M_2}\]
\[ \Updownarrow \]
\[ f_1(u,s) = 0.\]
But since $f_1$ is also a correction function for the $M_1$-choice
$\chxyz_1$,
\[ f_1(u,s) = 0\]
\[ \Updownarrow \]
\[ \langle a_0\dots a_k u_0\dots u_k x^1_{u_0s} \dots x^1_{u_{k-1}s}
y^1_{u_ks} z^1_u \rangle \in Q_s^{M_1}.\]
So, we have both $\bbz \models \sum_{\ell <k}i_\ell = h_k(u_0)+g(s)$
and $\bbz \models \sum_{\ell <k}i^1_\ell = h^1_k(u_0)+g^1(s)$, so
setting

\[ g_1(u_\ell,s) = i^1_\ell,\quad g_2(u_k,s) = h^1_k,\quad g_3(u) =
g^1\]
yields
\[ \bbz_2 \models \sum_{\ell<k} i_\ell - h_k(u_0) - g(s) = \sum
_{\ell<k} g_1(u_\ell,s) - g_2(u_k,s)(u_0) - g_3(u)(s). \]
Since $f_1$ does this for all possible $(k+1)$-tuples, we have
all the compability we need.

\item[(c) $\to$ (a)] If the predicates are the same modulo $g_1$, $g_2$
and $g_3$ then obtaining {\bf (a)} becomes a matter of building $F:M_1
\stackrel{\approx}{\longrightarrow}_{P_0^{M_1}\cup P_1^{M_1}}
M_2$. Clearly we can start by
$F\restriction P_0^{M_1} = id$, and then extend its definition to all
the other portions of the model. The only strong restriction to the
extension of this
to the whole model is given by the relations $Q_s^{M_1}$ and
$Q_s^{M_2}$. But part (4) of \textbf{(c)} provides this: the functions $g_1,g_2,g_3$ provide the
definition of the isomorphism. Precisely, let $\langle a_0\dots a_ku_0\dots u_kx_0\dots
x_{k-1}y_k,z\rangle$ be a tuple from $M_1$; we use (4) to find simultaneously $F(x_\ell)$, $F(y_k)$
and $F(z)$. Compute (in $\bbz_2$) the value $\sum_{\ell<k} i_\ell - h_k(u_0) - g(s)$ corresponding
to the tuple. For every $s\in S$, this value is $0$ iff the tuple belongs to $Q_s$. By (4), this
value is equal to $\sum_{\ell<k} g_1(u_\ell,s) - g_2(u_k,s)(u_0) - g_3(u)(s)$. But also by (4), this
value also corresponds to a corresponding tuple $\langle a_0\dots a_ku_0\dots u_kx_0'\dots
x_{k-1}'y_k',z'\rangle$ in $M_2$. This provides the values $F(x_\ell)=x_\ell'$, $F(y_k)=y_k'$ and
$F(z)=z'$: $\langle a_0\dots a_ku_0\dots u_kx_0\dots x_{k-1}y_k,z\rangle\in Q_s^{M_1}$ if
and only if $\langle a_0\dots a_ku_0\dots u_kx_0'\dots x_{k-1}'y_k',z'\rangle\in Q_s^{M_2}$.
\end{description} 
\qed

\begin{remark}
Counting the number of isomorphism types here is akin to the study of
$Ext(G,\bbz)$ in the work of the first author and V\"ais\"anen in
\cite{ShVa:646}\footnote{Here $I(\lambda,\psi)$ is counted by the group of
correction functions, derived from {\it some} $g_1$, $g_2$, $g_3$:
\[ I(\lambda,\psi_k^\lambda) = \Big \{ f \mbox{ a correction function} \Big|
  f(u,s) = \sum_{\ell<k}g_1(u_\ell,s) -g_2(u_0,s) -g_3(u)\Big \} .\]}.
\end{remark}

\begin{lemma}\label{modset2}
If $M_{I_1,f_1}$ and $M_{I_2,f_2}$ are models of $\psi$, and $h:I_1\to
I_2$ is one-to-one and onto, {\sf then} the following are equivalent:
\begin{description}
    \item[(a)] there is an isomorphism from $M_{I_1,f_1}$ onto $M_{I_2,f_2}$ extending $h$.
    \item[(b)$_1$] there is an $M_{I_2,f_2}$-choice $\chxyz$
whose correction function is $f_1$,
\item[(b)$_2$] there is an $M_{I_1,f_1}$-choice $\chxyz$
whose correction function is $f_2$,
\item[(c)] there are functions $g_1$, $g_2$, $g_3$ (``to
correct the choice of zeros''), with
\begin{enumerate}
\item $g_1:\ialak \times S \to \bbz_2$ (like the $x_{u,s}$'s above),
\item $g_2:\ialak \times S \to H_I$ (like the $y_{u,s}$'s above),
\item $g_3:\ialakmasuno  \to G$ (like the $z_u$'s above),
\item if $\langle a_0\dots a_ku_0\dots u_kx_0\dots
    x_{k-1}y_k,z\rangle$ are like in Definition \ref{defM_I} for $M_{I_1,f_1}$,
    or $M_{I_2,f_2}$ \underline{then} 
\[ \bbz_2 \models \sum_{\ell<k} i_\ell - h_k(u_0) - g(s) = \sum
_{\ell<k} g_1(u_\ell,s) - g_2(u_k,s)(u_0) - g_3(u)(s) \]
\end{enumerate}
\end{description}

%
%
\end{lemma}

\Proof The proof is almost the same as that of the previous lemma~(\ref{modset}). The main change is
that now the identity on $I$ is replaced by a bijection from $I_1$ onto $I_2$; the rest of the proof
amounts to a renaming via the bijection $F\restriction I_1$. \qed

An important special case of the previous lemma happens when $I_1=I=I_2$ \emph{but the isomorphism
is not the identity on $I$}. In this case, the restriction of the isomorphism between $M_{I,f_1}$
and $M_{I,f_2}$ is a \emph{permutation} of $I$. Our combinatorial criterion for non-isomorphism will
focus on this case.

Recall $\mathfrak D$
is the regular filter on $S$ generated by sets of the form
$\langle u\rangle = \{ v\in S | u\subset v\}$, where
$S=[\lambda]^{<\aleph_0}$: $$\gd = \gd_\l := \{ A\subset S | \exists
u_A\in S \forall v\in S(u_A\subset v \to v\in A)\}$$
  (see definition~\ref{bricks}).

  The notion of an $I$-function, which we define next, is central to our combinatorial criterion.

\begin{definition}
$f:[I]^{k+1}\times S\to \bbz_2$ is an $I$-function iff 
\[ \{ s\in S | f_u(s) \not= 0\} \in {\frak D}, \mbox{ for all }u\in
 [I]^{k+1},\]
where $f_u: S\to \bbz_2$ is given by $f_u(s) = f(u,s)$.\\
\end{definition}

\begin{lemma}\label{7.2}
    Let $f:[I]^{k+1}\times S\to \bbz_2$ be an $I$-function.
{\sf Then,} the following is a sufficient condition for
\[M_{I,f}\not\approx M_I:\]


\begin{description}
    \item[($\star$)] for every $F_1:\ialak \to [I]^{\leq \l}$, $F_2: \ialak \to
{}^S(\bbz_2)$ and $\pi$ a permutation of $I$, there exists 
$u=\{ t_0,\dots,t_k\}\in \ialakmasuno$ (i.e., with no repetitions) such that
\begin{description}
    \item[($\alpha$)] $t_k \not\in F_1(\{ t_0\dots t_{k-1}\})$,
\item[($\beta$)] $f_{\pi\{ t_0,\dots,t_k\}} - \sum_{\ell<k}
    F_2(\{ t_0,\dots,t_k\}\setminus \{ t_\ell\}) \not\in G$.
\end{description}
\end{description}
\end{lemma}

\noindent
Before proving Lemma~\ref{7.2}, we note:
\begin{itemize}
    \item First, $(\star)$ is a purely combinatorial statement; this will enable us to focus
	solely on combinatorial principles to prove the failure of categoricity.
    \item Also, by the definition of $G$, $(\beta)$ says that
	for $\mathfrak D$-few elements $s\in S$ do we have 
\[ f_{\pi(t_0,\dots,t_k)}(s) = \sum_{\ell<k}F_2(\{ t_0,\dots,t_k\}
\setminus \{ t_\ell\} )(s).\]
\end{itemize}

%
%
%
%
%


\noindent
Notice the role of the permutation $\pi$ of $I$ in the combinatorics
that follows.

\bigskip
\Proof of \ref{7.2}. Assume that $M_{I,f}\approx M_I$. Then, since $M_I\approx M_{I,{\mathfrak
0}}$ (the null correction function) we may apply Lemma~\ref{modset2} and (by (b)${}_2$ of that
lemma)
assume that $\chxyz$ witnesses $M_{I,f} \approx
M_{I,0}$, with correction function identically zero.

We construct $F_1,F_2$ such that $(\star)$ of~\ref{7.2} does not hold (for the permutation $\pi$
induced by the isomorphism between $M_{I,f}$ and $M_I$).

We first let $F_1: [I]^k \to [I]^{\leq \l}$ be
\[ F_1(v) = \bigcup \{ w\in \ialak | \mbox{ for some } s_1\in S,
y_{v,s_1}(w)\not= 0\} .\]
This is well defined, as $F_1(u)$ is a union of $|S|$-many  finite sets.
Also, let
\[ F_2(v) = \langle x_{v,s} | s\in S\rangle.\]

We will show that \textbf{no} $u\in \ialakmasuno$ satisfies both $(\alpha)$ and $(\beta)$ of
condition $(\star)$.

Suppose otherwise; let then $u=\left\{ t_0,\dots,t_k\right\} \in \ialakmasuno$
satisfy
$(\a)+(\b)$. Let as usual $u_\ell = u\setminus
\{ t_\ell\}$. By $(\a)$, for each $s$,
\[ y_{u_k,s} (u_0) = 0.\]

\noindent
[Why? Just notice that by $(\a)$, $$t_k\notin F_1(u_k)
= \bigcup \{ v\in \ialak | \mbox{ for some } s_1\in S,
y_{u,s_1}(v)\not= 0\},$$ so for all $v\in [I]^k$, if $t_k\in v$, then
for all $s_1\in S$ we have $y_{u_k,s_1}(v)=0$. In particular, as
$t_k\in u_0$, $y_{u_k,s_1}(u_0)=0$.]

Now, since $\chxyz$ is an $M_{I,f}$-choice with correction function
identically zero, for each $s\in S$ we have that
\[ \langle a_0,\dots,a_k,u_0,\dots,u_k,x_{u_0,s},\dots,x_{u_{k-1},s}, y_{u_k,s},z_u\rangle\]
belongs to $Q_s^{M_{I,f}}$ if and only if (by the definition of this predicate in the model
$M_{I,f}$)
\[ \bbz_2\models \sum_{\ell<k}x_{u_\ell,s}=y_{u_k,s}(u_0)+z_u(s)+f_{\pi(u)}(s).\]
But 


But we also have that $z_u(s) = 0$ for the $\gd$-majority of $s\in S$ (by the definition of
$G$); since we also have that $y_{u_k,s}(u_0) = 0$ for our particular $u$,
\begin{description}
\item[(*)] For the $\gd$-majority of $s\in S$
\[ \sum_{\ell<k} x_{u_\ell,s} = f_{\pi(u)}(s) .\]
\end{description}


But this contradicts $(\b)$. \qed

\begin{remark} 
We can then regard $F_2$ as
\[F_2:\ialak \to {}^S(\bbz_2)/G .\]
\end{remark}

\begin{corollary}
If $f_1$, $f_2$ are $I$-functions, and $f= f_1-f_2$ (coordinatewise)
satisfies $(\star)$, then $M_{I,f_1}\not\approx M_{I,f_2}$.
\end{corollary}

\Proof We already know that since $f$ satisfies $(\star)$, $M_{I,f}\not\approx M_I$. Now suppose
we have an isomorphism $F:M_{I,f_1}\stackrel{\approx}{\longrightarrow} M_{I,f_2}$. As before,
$\pi=F\restriction I$ is a permutation of $I$, and the automorphism lifts in a natural way to all
components of $M_{I,f}$ in the vocabulary $\tau^-$. Now, the remaining part of $\tau$: if $s\in
S$, then a tuple\\ $\langle
a_0,\dots,a_k,u_0,\dots,u_k,x_{u_0,s},\dots,x_{u_{k-1},s},y_{u_k,s},z_u\rangle$ belongs to $Q_s$
in $M_{I,f}$ if and only if (for the corresponding indices)
\[\bbz_2\models \sum_{\ell<k}i_\ell=h_k(u_0)+g(s)+f(u,s)\]
but this holds if and only if
\[\bbz_2\models \sum_{\ell<k}i_\ell=h_k(u_0)+g(s)+f_1(u,s)-f_2(u,s).\]
Now, since $F$ is an isomorphism, this is true if and only the tuple\\
$\langle
F(a_0),\dots,F(a_k),F(u_0),\dots,F(u_k),F(x_{u_0,s}),\dots,F(x_{u_{k-1},s}),
F(y_{u_k,s}),F(z_u)\rangle$
belongs to $Q_s$ in $M_I$, this is if and only if
\[\bbz_2\models \sum_{\ell<k}i_\ell'=h_k'(u_0')+g'(s)\]
where the primes denote the values corresponding to the $F$-images of components of the long
tuple. But this witnesses that $F$ is also an
isomorphism between $M_{I,f_1}$ and $M_{I,f_2}$, which contradicts the hypothesis. \qed

\medskip
We now have what we need for a proof of failure of categoricity at some $\mu$ above the
categoricity cardinals. Notice we do \emph{not} give an optimal (minimal) such $\mu$; this is left
for (possible) later work.

\begin{theorem}
    For some $\mu > \lambda^{+k}$, the sentence $\psi_k^\lambda$ is not categorical in $\mu$.
    \label{failcat}
\end{theorem}

\Proof Let $\mu$ be a cardinal with the following properties:

\begin{description}
\item[$\otimes_1$] $\mu \to (\w)^k_{2^\l}$,
\item[$\otimes_2$] $\mu \not\to (\w)^{k+1}_{2^\l}$,
\item[$\otimes_3$] $\mu$ regular.
\end{description}

The existence of such a $\mu$ uses the Erdös-Rado theorem (the partition
$$\beth_k(\lambda)^+\to \left( (2^\lambda)^+\right)^k_{2^\lambda}$$
is an instance) for $\otimes_1$ and the negative partition relation $\beth_{k+1}(\lambda)\not\to
(k+2)^{k+1}_{2^\lambda}$ (a consequence of~\cite[Lemma 24.1(e)]{EHMR}) for $\otimes_2$; we may
therefore take $\mu$ as $\beth_k(\lambda)^+$.


Let then $I$ have cardinality $\mu$ and let $f:[\mu]^{k+1}\to G^+/G$ be an $I$-function witnessing
$\otimes_2$ (recall that $G^+$ denotes the group ${}^S\bbz_2$). We use our criterion~\ref{7.2} to
show that $M_{I,f}$ and $M_I$ can not be isomorphic from which we conclude that the sentence
$\psi_k^\lambda$ is not categorical in $\mu$.

Let $F_1:\ialak \to [I]^{\leq \lambda}$, $F_2:\ialak \to G^+$ and $\pi$ a permutation of $I$. 

Now find $E\subset \mu$ club such that
\[\a_0<\dots <\a_k\in E \Longrightarrow 
\left \{ \begin{tabular}{l}
$F_1(\a_0,\dots ,\a_{k-1})\subset \a_k$,\\
$\pi(\a_0),\dots,\pi(\a_{k-1}) < \a_k$.
\end{tabular}
\right . \]

This is possible by the regularity of $\mu$.

Now apply $\otimes_1$ to $F_2\restriction E$: since $\mu\to(\omega)^k_{2^\lambda}$, there must be
an infinite $\omega$-sequence $T=\left\{ \alpha_0<\alpha_1<\dots \alpha_n<\dots\right\}$ such that
$F_2\restriction [T]^k$ is constant. Therefore we have, for $u=\left\{ \alpha_0,\dots,\alpha_k
\right\}$ and $u_\ell=u\setminus \left\{ \alpha_\ell \right\}$:
\begin{itemize}
    \item $\alpha_k\notin F_1(\left\{ \alpha_0,\dots,\alpha_{k-1} \right\})$ (since these are
	elements from the club $E$) and
    \item the equation $f_{\pi\{\alpha_0,\dots,\alpha_k\}}(s) = \sum_{\ell<k}F_2(u_\ell)(s)$ holds
	for $\gd$-few elements $s$: as $F_2$ is constant on $u_\ell$ from the monochromatic
	sequence, the sum on the right hand side will be $0$ when $k$ is even (and $1$ when $k$ is
	odd) whereas the value on the left hand side will not be constant, by $\otimes_2$ applied
	to $f$.
\end{itemize}
The previous two properties correspond to $(\alpha)$ and $(\beta)$ of the criterion from
Lemma~\ref{7.2}. Therefore, $M_{I,f}\not\approx M_I$ and the sentence $\psi^\lambda_k$ is not
categorical in $\mu$.
\qed

The result also holds for all $\kappa\geq \beth_{k+1}(\lambda)^+$ (we will show monotonicity of
the crucial criterion).

\begin{conclusion}
    The sentence $\psi^\lambda_k$ is not categorical in any $\kappa\geq\beth_{k+1}(\lambda)^+$.
\end{conclusion}

\Proof Let $\kappa \geq \mu= \beth_{k+1}(\lambda)^+$. If  $\kappa =
\mu$, Theorem~\ref{failcat} shows how to get \textbf{two} non-isomorphic models. If
$\kappa > \mu$ then let $J$ be a set of cardinality $\kappa$.

We show that as $\kappa>\mu$ we may pick a $J$-function $f:[J]^{k+1}\times S\to \bbz_2$ satisfying
the criterion $(\star)$ of Lemma~\ref{7.2} (which will enable us to conclude that
$M_{J,f}\not\approx M_J$, and thus conclude failure of categoricity at $\kappa$).

Let first $F_1:[J]^k\to [J]^{\leq \lambda}$, $F_2:[J]^k\to {}^S\bbz_2$ and $\pi$ a permutation of $J$. Let $I\subset J$ with $|I|=\mu$, $I$ closed under $\pi$ and such that $F_1\restriction
\ialak:\ialak \to [I]^{\leq \lambda}$. [Such an $I$ exists by closing first under iterating taking
the unions of $F_1$-images of $k$-tuples from $I$ and taking the union of $\mu$ many sets of
cardinality $\leq \lambda <\mu$ - after an $\omega-$iteration the result is closed under
$F_1$-images. Similarly, we close under images and preimages under the permutation $\pi$ and
alternate these closure operations $\omega$ many times.]

Let now $f:[J]^{k+1}\times S\to \bbz_2$ be a $J$-function that witnesses $\otimes_2$ on the set $I$;
as $|I|=\mu$, this is possible.

Furthermore, for the set $I$, the functions $F_1\restriction I$, $F_2$ and $\pi\restriction I$ are in the
situation of Lemma~\ref{7.2}. The proof of Theorem~\ref{failcat} applies then, as $|I|=\mu$. We
then obtain $u=\left\{ t_0,\dots,t_k \right\}\in \ialakmasuno$ such that $(\alpha)$ and
$(\beta)$ of the criterion hold. But these properties are also true of the original
$F_1,F_2,\pi$. Therefore, $M_{J,f}\not\approx M_J$. \qed

\begin{remark}
Here are some important differences between the structure of this
proof and that of \cite{HaSh:323}:
\begin{enumerate}
\item The use of the filter $\gd$ is central here - it is not needed there.
\item The way the group itself is used is slightly different at the
end of the proof.
\end{enumerate}
\end{remark}

We conjecture that the class has the maximal number of models at all
$\mu>\beth_{k+1}(\lambda)^+$.

\section{Further directions}

After our generalization of the original Hart-Shelah example to the stronger
logic $L_{(2^\lambda)^+,\omega}$, we have the following situation:
\begin{itemize}
	\item Any generalization of the early results from 1983 for the logic
		$L_{(2^\lambda)^+,\omega}$ must
		necessarily use as hypothesis few models in all cardinalities
		$\lambda,\lambda^+,\dots,\lambda^{+k},\dots$ for all
		$k<\omega$. The first author has written several papers in this
		direction (see~\cite{Sh09a}), in the (wider) context of AECs.
	\item On the other hand, the necessity of an interval of
		$\aleph_0$-many cardinals with few models to start the
		machinery for categoricity transfer seems interesting per se;
		and even more so the fact that this would happen along all
		strengthenings of $L_{\omega_1,\omega}$ (inside
		$L_{\infty,\omega}$).
	\item Finally, we conjecture that our sentence may be analyzed in terms
		of building frames, in the spirit of the work of Boney and
		Vasey~\cite{BoVa}. Specifically, we conjecture that our
		abstract elementary class
		\[ {\mathcal
    K}^*(\lambda,k)=(Mod(\psi_k^\lambda),\prec_{(2^\lambda)^+,\omega})\]
    is $(<\lambda,\lambda^{+k-1})$-tame,
  $(<\lambda,\lambda^{+k-1})$-typeshort over models of size
  $\lambda^{+k-2}$, and that
  \begin{enumerate}
\item for each $m\leq k-1$ there is a frame ${\mathfrak
    s}^*(\lambda,k)_m$ that is type-full and $\lambda^{+m}$-good on
  $Mod(\psi_k^\lambda)$,
\item The (type-full and $\lambda^{+k-1}$-good) frame ${\mathfrak
    s}^*(\lambda,k)_{k-1}$ is not weakly successful.
  \end{enumerate}
\end{itemize}

\bibliographystyle{elsarticle-num-names}

\bibliography{listb,listx}
\end{document}